\documentclass[11pt,a4paper]{amsart}
\usepackage[utf8]{inputenc}
\usepackage{graphicx}
\usepackage{amsmath,amssymb,amsfonts,amsthm}
\usepackage{booktabs}
\usepackage{tikz}
\usetikzlibrary{shapes.geometric}
\usepackage[hidelinks]{hyperref}
\theoremstyle{plain}
\newtheorem{theorem}{Theorem}[section]
\newtheorem{proposition}[theorem]{Proposition}
\newtheorem{lemma}[theorem]{Lemma}
\newtheorem{corollary}[theorem]{Corollary}
\newtheorem{maintheorem}{Theorem}

\theoremstyle{remark}
\newtheorem{remark}[theorem]{Remark}
\newtheorem{example}[theorem]{Example}
\theoremstyle{definition}
\newtheorem{definition}[theorem]{Definition}

\newcommand{\ZZ}{\mathbb Z}
\newcommand{\FF}{\mathbb F}
\newcommand{\PP}{\mathbb P}

\newcommand{\ima}{\mathbf i}
\DeclareMathOperator{\End}{End}

\numberwithin{equation}{section}
\AtBeginDocument{\interdisplaylinepenalty=10000}
\begin{document}
\title[Supersingular Legendre curves and isogeny cycles]{A recurrence characterization of supersingular Legendre curves
and a product--trace identity for 2-isogeny cycles}
\author{Juan Jes\'us Le\'on}
\address{GMV, C/ Isaac Newton 11, P.T.M., Tres Cantos, 28760, Madrid, Spain}
\email{jjleon@gmv.com}
\thanks{Juan Jes\'us Le\'on is the corresponding author. These authors contributed equally to this work.}
\author{Vicente Mu\~noz}
\address{Departamento de \'Algebra, Geometr\'{\i}a y Topolog\'{\i}a,
Universidad Complutense de Madrid, Ciudad Universitaria, 28040, Madrid, Spain}
\email{vicente.munoz@ucm.es}
\subjclass[2020]{11G20, 14H52, 11B37}
\keywords{Elliptic curves, isogeny graphs, finite fields, nonlinear recurrence,
arithmetic-geometric mean, supersingularity}
\begin{abstract}
We study the recurrence $v_{n+1}^2=(v_n+v_n^{-1})/2$ over
quadratic finite fields and its relation to explicit $2$-isogenies.
Starting from a root form of elliptic curves, we derive the recurrence
and relate it to the arithmetic-geometric mean. We prove that every
branch from a nonsingular initial parameter remains in the quadratic
field if and only if its associated Legendre curve is supersingular.
We also establish a finite criterion in terms of three Frobenius-stable
cyclic subgroups. For closed recurrence sequences, a normalized product
of the parameters is the differential multiplier of the composed
endomorphism. A transported point of order four fixes its trace modulo
four and yields an explicit product--trace identity. We discuss the
edge parameter of the $2$-isogeny graph, exceptional configurations,
and the distinction between this polynomial identity and the counting
of primitive cycles.
\end{abstract}
\maketitle

%%%%%%%%%%%%%%%%%%%%%%%%%%%%%
\section{Introduction}\label{sec:intro}
%%%%%%%%%%%%%%%%%%%%%%%%%%%%%
Isogenies relate the geometry of elliptic curves to the arithmetic of
their endomorphisms. Fixing their degree produces a graph in which
vertices represent geometric isomorphism classes and edges arise from
quotients by finite subgroups. In degree two, an elliptic curve over
an algebraically closed field of characteristic different from two
has three possible kernels. After one isogeny has been chosen, two
choices remain if the next step is required not to be its dual.
Explicit equations for these choices provide a way to describe paths
by algebraic recurrences. Such a description raises two questions:
over which field can the recurrence be continued, and what arithmetic
information is retained when a sequence closes up to a cycle?

In this paper we study these questions for the recurrence
\begin{equation}\label{eq:introrecurrence}
 v_{n+1}^2=\frac12\left(v_n+\frac1{v_n}\right)
\end{equation}
in characteristic $p\geq5$. We write $k=\FF_{p^2}$ and
$\bar k=\overline{\FF}_p$. A parameter is called admissible when
$v\ne0$ and $v^4\ne1$. Each admissible parameter has two distinct
successors over $\bar k$, obtained by choosing either square root,
and both successors are again admissible. Thus a starting value
determines a family of branches. Even if it belongs to $k$, it is
not automatic that all subsequent square roots belong to $k$.
Our first main result characterizes exactly when every branch
remains in this quadratic field. Our second main result identifies
the product along a closed branch with a scalar arising from an
elliptic-curve endomorphism.

The geometric interpretation uses the Legendre models
\[
 E_v:\quad y^2=x(x-1)(x-v^4),\qquad
 J(v)=j(E_v)=256\frac{(v^8-v^4+1)^3}{(v^8-v^4)^2}.
\]
For every step $v\mapsto w$ of~\eqref{eq:introrecurrence}, we
give an explicit $2$-isogeny $\Psi_v:E_v\to E_w$ with kernel
$\{O,(0,0)\}$. Successive maps do not backtrack. The parameter
$v$ carries more information than the invariant $J(v)$: for
example, $J(v)=J(-v)$, although signs affect subsequent recurrence
steps. Consequently, paths in the parameter graph must be
distinguished from their projections to the graph of $j$-invariants.
We also keep track of the chosen curve models and maps, since equality
of invariants gives an isomorphism over $\bar k$ but need not give
one over the field of definition; see~\cite{Silverman}.

We obtain the recurrence initially from a different description,
the root form of an elliptic curve. Starting with the quotient by
a specified point of order two, we normalize its equation and
write the two forward kernel coordinates explicitly. A sequence
of changes of variables then leads to~\eqref{eq:introrecurrence}.
The quotient calculations use the explicit theory of
$2$-isogenies~\cite{Vel,Shu}; we include the intermediate
substitutions to make the relation between kernel coordinates,
recurrence parameters and Legendre equations explicit. The root
normalization has exceptional values, whereas the Legendre
description applies to every admissible parameter. The congruence
condition $p\equiv7\pmod{12}$ used in the root discussion and
some examples is therefore not a hypothesis of the two main
theorems.

The recurrence also has a direct interpretation in terms of the
arithmetic-geometric mean. For a branch $(v_n)$, put
$b_0=1$, $b_{n+1}=v_nb_n$ and $a_n=v_n^2b_n$. Then
\[
 a_{n+1}=\frac{a_n+b_n}{2},\qquad b_{n+1}^2=a_nb_n,
 \qquad v_n=\frac{b_{n+1}}{b_n}.
\]
Thus the signs of the square roots in the recurrence are precisely
the choices of geometric mean in the associated AGM sequence.
Over the positive reals, choosing positive roots gives the
classical convergent process studied in~\cite{Cox}; over a finite
field, retaining both choices gives a directed graph.
Griffin, Ono, Saikia and Tsai~\cite{REF-1} study finite-field AGM
graphs and their elliptic curves, and Kayath, Lane, Neifeld, Ni and
Xue~\cite{AGMaquariums} develop this correspondence over arbitrary
finite fields. Our use of the AGM makes the field-of-definition
problem for~\eqref{eq:introrecurrence} accessible through their
results on the structure of these graphs. The real seed
$(a_0,b_0)=(\sqrt2,1)$ corresponds to $v_0=\sqrt[4]{2}$ and
$J(v_0)=1728$; see Example~\ref{ex:seed}.

The first main statement is the following characterization,
proved as Theorem~\ref{thm:A}.
\medskip

\noindent\textbf{Theorem A.}
\emph{Let $v_0\in k$ be admissible. Every finite branch of
\eqref{eq:introrecurrence} starting at $v_0$ lies in $k$ if and
only if $E_{v_0}$ is supersingular.}
\medskip

The quantifier over all branches is essential: the assertion
concerns every successive choice of signs, not a single selected
sequence. Supersingular invariants are known to belong to
$\FF_{p^2}$, but this alone does not establish rationality of the
parameters and roots used in a recurrence. The proof uses the
AGM correspondence explicitly. In the supersingular case the
required successors remain rational. Conversely, if all branches
remain in $k$, they produce a finite nonempty set of AGM vertices
closed under both successors. The structure theorem for ordinary
AGM components in~\cite{AGMaquariums} rules out such a set in the
ordinary case. We also give a separate finite criterion using
three Frobenius-stable cyclic subgroups of $E[2^n]$ in a specified
relative position. This criterion requires only sufficiently long
finite chains, but it does not assert that three arbitrary chains
suffice.

For the second main result, the normalization of the maps is
important. With $\xi_v=dx/(2y)$, their pullback relation is
\[
 \Psi_v^*\xi_w=2v\xi_v.
\]
An admissible recurrence cycle of length $N$ is a sequence
$(v_0,\ldots,v_{N-1})$ satisfying the recurrence with indices
read modulo $N$. It returns to the same parameter and hence to
the same chosen model. Its composite map and parameter product are
\[
 \Theta=\Psi_{v_{N-1}}\circ\cdots\circ\Psi_{v_0}
 \in\End(E_{v_0}),\qquad
 \omega=\prod_{n=0}^{N-1}v_n.
\]
The differential multiplier of $\Theta$ is therefore $2^N\omega$.
This observation connects the recurrence product with the
characteristic polynomial of an endomorphism. Theorem~\ref{thm:B}
gives the following precise statement.
\medskip

\noindent\textbf{Theorem B.}
\emph{For an admissible recurrence cycle over $\bar k$ of length
$N\geq2$, the endomorphism $\Theta$ has degree $2^N$ and integer
trace $t$ satisfying}
\[
 t\equiv1\pmod4,\qquad |t|\leq2^{1+N/2},\qquad
 2^N\omega^2-t\omega+1=0.
\]
\emph{The last identity is interpreted in characteristic $p$,
with $t$ reduced modulo $p$; its roots belong to $\FF_{p^2}$.}
\medskip

The trace congruence comes from an explicit point of order four
transported by each map and fixed by the composite around the
cycle. Applying the characteristic equation of $\Theta$ to this
point gives $t\equiv1\pmod4$, while its action on differentials
gives the quadratic identity. Together with the trace bound,
this yields a polynomial, depending only on $N$, satisfied by
every cycle product. It constrains the possible products but
does not count primitive cycles: different cycles can share a
product, and a closed sequence of length $N$ can have smaller
least period. Reversing the sequence and applying
$T(v)=(v+1)/(v-1)$ produces another cycle with product
$1/(2^N\omega)$, the other root of the quadratic, counted with
multiplicity. Frobenius conjugation instead gives $\omega^p$;
these two operations have the same effect on the product when
$\omega\notin\FF_p$.

The passage from parameters to geometric invariants also motivates
our treatment of exceptional graph configurations. Using the
classical modular polynomial $\Phi_2$, we introduce an edge label
$g$ for which adjacent invariants $X,Y$ satisfy
$XY=g^3$ and $X+Y=g^2-495g+54000$. These relations give explicit
descriptions of loops and of distinct edge labels with the same
endpoints. The diagrams distinguish labels from kernel
multiplicities, which can behave differently at exceptional
vertices; related analyses appear in~\cite{Ghantous,REF-5}.
This distinction is also relevant to cycles, for which the
specified maps retain information lost on passing to invariants.

The paper is organized as follows. Section~\ref{sec:prelim}
recalls the elliptic-curve notation and the explicit quotient and
dual $2$-isogenies. Section~\ref{sec:root} introduces the root
form and derives the recurrence through successive changes of
variables. Section~\ref{sec:legendre} constructs the Legendre
models and proves the isogeny and differential formulas.
Section~\ref{sec:agm} develops the correspondence with the AGM.
Section~\ref{sec:A} proves Theorem A and the finite
criterion involving three Frobenius-stable directions.
Section~\ref{sec:graph} studies the modular equation, the edge
parameter and the exceptional configurations.
Finally, Section~\ref{sec:B} proves Theorem B, derives the
cycle-product polynomial, and examines reversal and Frobenius
conjugation.

%%%%%%%%%%%%%%%%%%%%%%%%%%%%%
\section{Elliptic curves and explicit 2-isogenies}\label{sec:prelim}
%%%%%%%%%%%%%%%%%%%%%%%%%%%%%
Throughout, $p\geq5$ is prime, $k=\FF_{p^2}$ and
$\bar k=\overline{\FF}_p$. For the root-form discussion and the
supersingular examples we assume $p\equiv7\pmod {12}$.
Thus $k=\FF_p[\ima]$, with $\ima^2=-1$, and $j=0$ is not supersingular.
Theorems~\ref{thm:A} and~\ref{thm:B}, formulated in Legendre coordinates,
do not require this congruence restriction.

An elliptic curve in short Weierstrass form is
\[
 E:y^2=x^3+ax+b,\qquad
 \Delta(E)=-16(4a^3+27b^2)\ne0.
\]
Its invariant is
\begin{equation}\label{eq:jweier}
 j(E)=1728\frac{4a^3}{4a^3+27b^2}.
\end{equation}
Equality of $j$-invariants characterizes isomorphism over $\bar k$,
not necessarily over $k$. We distinguish a chosen model over $k$ from
its geometric isomorphism class. We use the standard facts that
supersingularity is invariant under isogeny and that all supersingular
$j$-invariants belong to $k$; see~\cite{Silverman}.

Every $2$-isogeny considered here is automatically separable.
Indeed, in characteristic $p$ the inseparable degree of an isogeny
is a power of $p$ dividing its total degree. Since $p\geq5$,
an isogeny of degree $2$ has inseparable degree $1$. We shall
therefore simply refer to these maps as $2$-isogenies.

\begin{proposition}\label{prop:quotient}
Let $P=(r,0)$ be a nonzero point of order two on $E$.
The quotient by $\langle P\rangle$, followed by a change of scale, has
$x$-coordinate
\begin{equation}\label{eq:quotient}
 \bar\psi(x)=\alpha\left(x+\frac{3r^2+a}{x-r}\right),\qquad \alpha\ne0,
\end{equation}
and codomain 
\begin{equation}\label{eq:quotientcurve}
 E' : y^2=x^3+\alpha^2(-4a-15r^2)x+
                 \alpha^3(-8ar-22r^3).
\end{equation}
For $c^2=\alpha$, the corresponding $y$-coordinate is
\[
 y'=c^3 y\left(1-\frac{3r^2+a}{(x-r)^2}\right).
\]
These formulas define a $2$-isogeny with kernel $\langle P\rangle$,
 $$
 \psi: E \longrightarrow E'
$$
and it is defined over any field containing the
coefficients, $r$ and $c$.
\end{proposition}

\begin{proof}
For $\alpha=1$, these are V\'elu's formulas for the kernel
$\langle P\rangle=\{O,P\}$; see~\cite{Vel,Shu}. The general case follows by the change
$(x,y)\mapsto(c^2x,c^3y)$ on the quotient. In particular, the scale
on $y$ must be specified when discussing the field of definition.
\end{proof}

The dual kernel has $x$-coordinate $-2\alpha r$. Indeed, the two
points of $E[2]\setminus\langle P\rangle$ have this common image.
The map induced by the dual on $x$-coordinates is
\begin{equation}\label{eq:dual}
 \bar{\hat\psi}(x)=\frac{x}{4\alpha}
 -\frac{3\alpha r^2/4+a\alpha}{x+2\alpha r}.
\end{equation}
Its composition with~\eqref{eq:quotient} is
\[
 \overline{[2]}(x)=\frac{x^4-2ax^2-8bx+a^2}{4(x^3+ax+b)}.
\]
A chain is called non-backtracking if each next kernel differs from
the kernel of the dual isogeny. This is a condition on
kernels on the specified models, even when some $j$-invariants coincide.

%%%%%%%%%%%%%%%%%%%%%%%%%%%%%
\section{Root form and the recurrence}\label{sec:root}
%%%%%%%%%%%%%%%%%%%%%%%%%%%%%
\subsection{The root form}
%%%%%%%%%%%%%%%%%%%%%%%%%%%%%

Several familiar models parameterized by $j$ exclude $j=1728$;
see~\cite{Con,Has,Was}. To include this value, for $j\ne0$ we use
the \textbf{root form}
\begin{equation}\label{eq:root}
 E(j):y^2=x^3-\frac{j}{576}x
              -\frac{j\sqrt{j-1728}}{20736\sqrt3}.
\end{equation}
The roots are chosen in $\bar k$. The two choices of the sign of
the constant coefficient give geometrically isomorphic curves.
Substitution in~\eqref{eq:jweier} gives $j(E(j))=j$, and
\begin{equation}\label{eq:rootdiscriminant}
 4a^3+27b^2=-\frac{j^2}{27648},\qquad
 \Delta(E(j))=\frac{j^2}{1728}.
\end{equation}
In particular, $E(1728):y^2=x^3-3x$. A root form may require a
larger field than a different model with the same $j$-invariant.

\begin{proposition}\label{prop:rootisogeny}
Let $P=(r,0)\in E(j)[2]$, with $r\ne0$, and put $\epsilon=j/r^2$.
On the locus where the following expressions are defined and the
codomain has nonzero $j$-invariant, take
\begin{equation}\label{eq:alpha}
 \alpha=\frac{12\sqrt3(2160-\epsilon)}
 {r(432-\epsilon)\sqrt{1728-\epsilon}}.
\end{equation}
The codomain of~\eqref{eq:quotient} is in root form, with
\begin{align}
 j'&=1728\frac{(2160-\epsilon)^3}
 {(432-\epsilon)^2(1728-\epsilon)},\label{eq:rootj}\\
 \sqrt{j'-1728}&=\frac{2592(1584-\epsilon)}
 {(432-\epsilon)\sqrt{1728-\epsilon}}.\label{eq:rootjminus}
\end{align}
The square roots and the sign of the root model are chosen compatibly.
\end{proposition}

\begin{proof}
Using $a=-j/576$, substitution in~\eqref{eq:quotientcurve} gives
\begin{align*}
 a'&=-3\frac{(2160-\epsilon)^3}
 {(432-\epsilon)^2(1728-\epsilon)},\\
 b'&=\frac{216}{\sqrt3}
 \frac{(2160-\epsilon)^3(\epsilon-1584)}
 {(432-\epsilon)^3(1728-\epsilon)^{3/2}}.
\end{align*}
Thus $a'=\displaystyle \frac{-j'}{576}$ and
$b'=\displaystyle \frac{-j'\sqrt{j'-1728}}{20736\sqrt3}$ with the displayed values.
Consequently the codomain equation is
\[
 E':\quad y^2=x^3-\frac{j'}{576}x
 -\frac{j'\sqrt{j'-1728}}{20736\sqrt3},
\]
which is precisely the root form~\eqref{eq:root} with parameter $j'$.
\end{proof}

For $r=0$ one has $j=1728$. The quotient is
$y^2=x^3+12\alpha^2x$, which has the same root form as the source
when $\alpha=\pm\ima/2$. At $j'=0$ the root chart must instead be
replaced by another Weierstrass model. None of these chart choices
is needed in the Legendre proofs below.

\begin{proposition}\label{prop:continuations}
In the notation of Proposition~\ref{prop:rootisogeny}, put
$\eta=\tfrac1{12} \sqrt{1728-\epsilon}$. The three order-two coordinates on
the codomain are
\begin{equation}\label{eq:twopoints}
 -2\alpha r,\qquad \alpha r(1+\eta),\qquad \alpha r(1-\eta).
\end{equation}
The first is the dual kernel; the remaining two are the forward kernels.
\end{proposition}

\begin{proof}
The monic cubic with these roots is
\[
 x^3-\alpha^2r^2(3+\eta^2)x+2\alpha^3r^3(1-\eta^2).
\]
Substituting $\eta^2=(1728-\epsilon)/144$, this polynomial becomes
\[
 x^3+\alpha^2r^2\left(\frac{\epsilon}{144}-15\right)x
 +\alpha^3r^3\left(\frac{\epsilon}{72}-22\right).
\]
Since $\epsilon r^2=j$ and $a=-j/576$, its coefficients are
$\alpha^2(-4a-15r^2)$ and $\alpha^3(-8ar-22r^3)$,
exactly those in~\eqref{eq:quotientcurve}. The assertion about the dual follows
from~\eqref{eq:dual}.
\end{proof}

%%%%%%%%%%%%%%%%%%%%%%%%%%%%%
\subsection{Changes of variables}\label{sec:changes}
%%%%%%%%%%%%%%%%%%%%%%%%%%%%%
We first specify the chain and its indexing, and then derive the
recurrence one change of variable at a time.

\begin{definition}\label{def:rootchain}
A \emph{root isogeny} is a $2$-isogeny between two chosen
curves in root form~\eqref{eq:root}. A \emph{root isogeny chain} is
a finite or infinite sequence
\[
 E_0\xrightarrow{\psi_0}E_1\xrightarrow{\psi_1}E_2
 \xrightarrow{\psi_2}\cdots,\qquad E_n=E(j_n),
\]
in which every $\psi_n$ is a root isogeny and, for each consecutive
pair of maps, $\ker\psi_{n+1}\ne\ker\widehat\psi_n$ on $E_{n+1}$.
Thus the chain is non-backtracking in the sense of Section~\ref{sec:prelim}.
\end{definition}

Write $j_n=j(E_n)$ and let $P_n=(r_n,0)$ be the nonzero point of
order two generating $\ker\psi_n$. In particular, $r_n$ is a
coordinate on the source $E_n$, whereas $r_{n+1}$ is a coordinate
on the target $E_{n+1}$ and specifies the \emph{next} isogeny.
Let $\alpha_n$ be the scale in~\eqref{eq:quotient} for $\psi_n$,
chosen as in Proposition~\ref{prop:rootisogeny}. We recall the
auxiliary parameters
\begin{equation}\label{eq:chainparameters}
 \epsilon_n=\frac{j_n}{r_n^2},\qquad
 \eta_n=\frac{\sqrt{1728-\epsilon_n}}{12},\qquad
 \epsilon_n=1728-144\eta_n^2.
\end{equation}
Here the sign of $\eta_n$ orders the two forward kernels; it can
be chosen independently of the already fixed map $\psi_n$.

For these calculations we work over $\bar k$, on a segment where
$r_n\ne0$ and the root normalizations in
Proposition~\ref{prop:rootisogeny} are defined. In particular,
$\epsilon_n\notin\{432,1728,2160\}$, and we retain a forward
kernel only while its coordinate $r_{n+1}$ is nonzero.
The cases in which a root chart fails are treated as explained
in Section~\ref{sec:root}; the Legendre recurrence in the next
section will not require these root charts.

%%%%%%%%%%%%%%%%%%%%%%%%%%%%%
\paragraph*{The recurrence for $\epsilon_n$ and $\eta_n$.}
%%%%%%%%%%%%%%%%%%%%%%%%%%%%%
Proposition~\ref{prop:continuations}, applied to $\psi_n$, gives
the three nonzero points of order two on $E_{n+1}$:
\[
 Q_{n+1,0}=(-2\alpha_nr_n,0),\qquad
 Q_{n+1,\pm}=(\alpha_nr_n(1\pm\eta_n),0).
\]
The point $Q_{n+1,0}$ generates $\ker\widehat\psi_n$. Since the
chain does not backtrack, $P_{n+1}$ is one of $Q_{n+1,+}$ and
$Q_{n+1,-}$. Consequently,
\begin{equation}\label{eq:nextrootcoordinate}
 r_{n+1}=\alpha_nr_n(1\pm\eta_n),\qquad
 r_{n+1}^2=\alpha_n^2r_n^2(1\pm\eta_n)^2.
\end{equation}

To compute $\epsilon_{n+1}$, first square~\eqref{eq:alpha} and
multiply by $r_n^2$:
\begin{equation}\label{eq:alphasquaredchain}
 \alpha_n^2r_n^2
 =432\frac{(2160-\epsilon_n)^2}
 {(432-\epsilon_n)^2(1728-\epsilon_n)}.
\end{equation}
On the other hand, formula~\eqref{eq:rootj} for the target invariant is
\begin{equation}\label{eq:targetjscale}
 j_{n+1}=1728\frac{(2160-\epsilon_n)^3}
 {(432-\epsilon_n)^2(1728-\epsilon_n)}
=4\alpha_n^2r_n^2(2160-\epsilon_n).
\end{equation}
Equivalently,~\eqref{eq:nextrootcoordinate} can be written as
\[
 r_{n+1}^2=\frac{j_{n+1}}{4(2160-\epsilon_n)}
 (1\pm\eta_n)^2.
\]
Now the definition of $\epsilon_{n+1}$ uses both quantities on
$E_{n+1}$. Substituting and cancelling the common nonzero factor
$\alpha_n^2r_n^2$ gives
\begin{equation}\label{eq:epsilonrecurrence}
 \epsilon_{n+1}
 =\frac{j_{n+1}}{r_{n+1}^2}
 =\frac{4\alpha_n^2r_n^2(2160-\epsilon_n)}
 {\alpha_n^2r_n^2(1\pm\eta_n)^2}
 =\frac{4(2160-\epsilon_n)}{(1\pm\eta_n)^2}.
\end{equation}

Using~\eqref{eq:chainparameters} at both indices, this becomes
$$
 1728-144\eta_{n+1}^2
 =\frac{4(2160-1728+144\eta_n^2)}{(1\pm\eta_n)^2} =\frac{576(\eta_n^2+3)}{(1\pm\eta_n)^2}.
$$
After rearrangement, we find
$$
 \eta_{n+1}^2 =12-\frac{4(\eta_n^2+3)}{(1\pm\eta_n)^2}=\frac{12(1\pm2\eta_n+\eta_n^2)-4\eta_n^2-12}
 {(1\pm\eta_n)^2} =8\eta_n\frac{\eta_n\pm3}{(\eta_n\pm1)^2}.
$$
We can use the plus sign uniformly by choosing the sign of
$\eta_n$ to record the selected forward kernel. Indeed, if the
minus sign occurs, set $\widetilde\eta_n=-\eta_n$. Then
\[
 8\eta_n\frac{\eta_n-3}{(\eta_n-1)^2}
 =8\widetilde\eta_n
 \frac{\widetilde\eta_n+3}{(\widetilde\eta_n+1)^2}.
\]
This change interchanges the labels $Q_{n+1,+}$ and $Q_{n+1,-}$,
without changing the chosen kernel itself. It also preserves the
preceding recurrence equation, which involves $\eta_n$ only
through its square. With these sign choices we therefore have
\begin{equation}\label{eq:eta}
 \eta_{n+1}^2=8\eta_n\frac{\eta_n+3}{(\eta_n+1)^2}.
\end{equation}

The same parameter recovers the target invariant. Indeed,
$ 2160-\epsilon_n=144(\eta_n^2+3)$,
$432-\epsilon_n=144(\eta_n^2-9)$ and
$ 1728-\epsilon_n=144\eta_n^2$.
Substitution in~\eqref{eq:rootj} cancels the three factors and yields
\begin{equation}\label{eq:jeta}
 j_{n+1}=1728\frac{(\eta_n^2+3)^3}
 {\eta_n^2(\eta_n^2-9)^2}.
\end{equation}

%%%%%%%%%%%%%%%%%%%%%%%%%%%%%
\paragraph*{From $\eta_n$ to $\mu_n$.}
%%%%%%%%%%%%%%%%%%%%%%%%%%%%%

Introduce the change of variable
\begin{equation}\label{eq:etatomu}
 \mu_n=\frac{\eta_n+3}{2\eta_n}
 =\frac12+\frac3{2\eta_n},\qquad
 \eta_n=\frac3{2\mu_n-1}.
\end{equation}
To express~\eqref{eq:jeta} in this variable, we compute
\[
 \eta_n^2=\frac9{(2\mu_n-1)^2},\quad
 \eta_n^2+3=\frac{12(\mu_n^2-\mu_n+1)}{(2\mu_n-1)^2},\quad
 \eta_n^2-9=\frac{-36\mu_n(\mu_n-1)}{(2\mu_n-1)^2}.
\]
The factors $(2\mu_n-1)^6$ cancel, giving
\begin{equation}\label{eq:jmu}
 j_{n+1}
 =\frac{1728\cdot12^3}{9\cdot36^2}
 \frac{(\mu_n^2-\mu_n+1)^3}{\mu_n^2(\mu_n-1)^2}
 =256\frac{(\mu_n^2-\mu_n+1)^3}{(\mu_n^2-\mu_n)^2}.
\end{equation}

For the recurrence, observe that
\[
 \frac{\eta_n+3}{\eta_n}=2\mu_n,\qquad
 \frac{\eta_n+1}{\eta_n}
 =1+\frac{2\mu_n-1}{3}=\frac{2(\mu_n+1)}3.
\]
Hence~\eqref{eq:eta} takes the form
\[
 \eta_{n+1}^2
 =\frac{8(\eta_n+3)/\eta_n}{((\eta_n+1)/\eta_n)^2}
 =\frac{16\mu_n}{4(\mu_n+1)^2/9}
 =\frac{36\mu_n}{(\mu_n+1)^2}.
\]
Choosing a square root $\sqrt{\mu_n}$ compatibly with
$\eta_{n+1}$, we may thus write
$\eta_{n+1}=6\sqrt{\mu_n}/(\mu_n+1)$. Substituting this expression
in~\eqref{eq:etatomu} at index $n+1$ gives
\begin{equation}\label{eq:mu}
 \mu_{n+1}
 =\frac12+\frac3{2\eta_{n+1}}
 =\frac12+\frac{3(\mu_n+1)}{12\sqrt{\mu_n}}
 =\frac12+\frac{\mu_n+1}{4\sqrt{\mu_n}}.
\end{equation}
The two square roots give the two choices of sign for
$\eta_{n+1}$ and hence the two corresponding values of $\mu_{n+1}$.
On the root-chart segment under consideration $\mu_n+1\ne0$:
by~\eqref{eq:etatomu} this is equivalent to $\eta_n+1\ne0$,
which follows from $r_{n+1}\ne0$ after the above sign convention.

%%%%%%%%%%%%%%%%%%%%%%%%%%%%%
\paragraph*{From $\mu_n$ to $u_n$ and $v_n$.}
%%%%%%%%%%%%%%%%%%%%%%%%%%%%%

Write $\mu_n=u_n^2$, choosing $u_n$ to be the square root used
in~\eqref{eq:mu}. That equation now reads
\[
 u_{n+1}^2
 =\frac12+\frac{u_n^2+1}{4u_n}
 =\frac{2u_n+u_n^2+1}{4u_n}
 =\frac{(u_n+1)^2}{4u_n}.
\]
After choosing a square root $v_n$ of $u_n$ with the compatible
sign, this is equivalent to
$u_{n+1}=\displaystyle \frac{u_n+1}{2v_n}$. Replacing $u_n$ and $u_{n+1}$ by
$v_n^2$ and $v_{n+1}^2$, respectively, gives
\begin{equation}\label{eq:recurrence}
 \boxed{ v_{n+1}^2=\frac{v_n^2+1}{2v_n}
 =\frac12\left(v_n+\frac1{v_n}\right).}
\end{equation}
Thus the successive changes of variable have led us to
$$
 u_{n+1}=\frac{u_n+1}{2\sqrt{u_n}},\qquad
v_{n+1}^2=\frac12(v_n+v_n^{-1}),
$$
and the compatible roots can also be read from consecutive terms:
\[
 u_n=\frac{\mu_n+1}{4\mu_{n+1}-2},\qquad
 v_n=\frac{u_n+1}{2u_{n+1}},
\]
whenever these denominators are nonzero.
\begingroup
More precisely, three consecutive ordered Legendre parameters
determine the recurrence parameter by
\[
 v_n=\frac{(\mu_n+4\mu_{n+1}-1)(2\mu_{n+2}-1)}
 {2(2\mu_{n+1}-1)(\mu_{n+1}+1)},
\]
whenever the preceding recovery formulas are defined. Thus an
ordered cycle determines all its recurrence parameters under this
nonvanishing condition, with indices read modulo its length.
For a finite chain $\mu_0,\ldots,\mu_M$, however, this argument
recovers only $v_0,\ldots,v_{M-2}$, since two subsequent Legendre
parameters are needed. The terminal parameters need not be unique:
changing the sign of $v_M$ preserves both its fourth power and
the last recurrence equation.
\par\endgroup
All these calculations
are over $\bar k$. Rationality of every branch over $k$ is a
separate assertion, proved later in Theorem~\ref{thm:A}.

Finally, $\mu_n=u_n^2=v_n^4$, so~\eqref{eq:jmu} becomes
\[
 j_{n+1}=256\frac{(v_n^8-v_n^4+1)^3}{(v_n^8-v_n^4)^2}.
\]
It is therefore convenient to introduce
\begin{equation}\label{eq:J}
 J(v)=256\frac{(v^8-v^4+1)^3}{(v^8-v^4)^2}.
\end{equation}
The index shift is worth keeping explicit: $j_n$ denotes the
invariant of the source $E_n$, whereas $J(v_n)=j_{n+1}$ is the
invariant of its target in this root-chain indexing. From now on
we index the associated Legendre curves directly by $v_n$, so
that their invariants are $J(v_n)$.

\begin{remark}\label{rem:level}
The parameter retains more information than an immediate choice of
an edge at the projected $j$-vertex. In particular $J(v)=J(-v)$.
In the changes of variables above, changing the sign of $v_n$
changes the sign of $u_{n+1}$, then the sign of
$2\mu_{n+2}-1$, and hence the sign of $\eta_{n+2}$.
Thus the corresponding choice of forward root kernel is delayed
by two levels. One must not identify the two signs at a recurrence
step with two different immediate successor $j$-invariants.
\end{remark}

%%%%%%%%%%%%%%%%%%%%%%%%%%%%%
\section{The Legendre Model}\label{sec:legendre}
%%%%%%%%%%%%%%%%%%%%%%%%%%%%%

\begin{definition}\label{def:admissible}
Over a field of characteristic different from $2$, a recurrence
parameter is admissible if $v\ne0$ and $v^4\ne1$.
Set 
 $$
 E_v:y^2=x(x-1)(x-v^4).
$$ 
A branch from $v_0$ is any sequence
in $\bar k$ satisfying~\eqref{eq:recurrence}; choices of both roots
are allowed at every step.
\end{definition}

The equation defining $E_v$ is in \emph{Legendre form}. Its
$j$-invariant is $J(v)$ from~\eqref{eq:J}. On the root charts of
Section~\ref{sec:changes}, where $J(v)\ne0$, the corresponding
root model has equation
\begin{equation}\label{eq:rootmodelv}
 E(J(v)):\quad
 Y^2=X^3-\frac{J(v)}{576}X
 -\frac{J(v)\sqrt{J(v)-1728}}{20736\sqrt3}.
\end{equation}
Thus $E_v$ and $E(J(v))$ are geometrically isomorphic models of
the same elliptic curve. In the indexing of the preceding section,
\begin{equation}\label{eq:legendreindexing}
 j(E_{v_n})=J(v_n)=j_{n+1},\qquad
 E_{v_n}\cong_{\bar k}E_{n+1}.
\end{equation}
Accordingly, a recurrence step denoted by $v\mapsto w$ below
means $v=v_n$ and $w=v_{n+1}$, with
$w^2=(v+v^{-1})/2$. The map $E_v\to E_w$ therefore corresponds
to the step from $E_{n+1}$ to $E_{n+2}$ at the level of geometric
isomorphism classes.

We can make the change of model explicit, including what happens
to the chosen kernel. With the sign convention of
Section~\ref{sec:changes}, the three roots of the cubic defining
$E_{n+1}$ are
\[
 r_{n+1}=\alpha_nr_n(1+\eta_n),\qquad
 q_{n+1}=\alpha_nr_n(1-\eta_n),\qquad
 t_{n+1}=-2\alpha_nr_n.
\]
Here $(r_{n+1},0)=P_{n+1}$ is the chosen kernel for the next
root isogeny, and $(t_{n+1},0)$ generates the dual kernel of
the preceding one. Set
\[
 d_n=q_{n+1}-r_{n+1}=-2\alpha_nr_n\eta_n\ne0,
 \qquad \delta_n^2=d_n,
\]
choosing $\delta_n$ in $\bar k$. In the root coordinates $(X,Y)$,
the equation of $E_{n+1}$ is
\[
 Y^2=(X-r_{n+1})(X-q_{n+1})(X-t_{n+1}).
\]
The affine change of coordinates
\begin{equation}\label{eq:roottolegendre}
 x=\frac{X-r_{n+1}}{d_n},\qquad
 y=\frac{Y}{\delta_n^3}
\end{equation}
transforms this equation into
\[
 y^2=x(x-1)(x-\mu_n),
\]
where
\[
 \mu_n=\frac{t_{n+1}-r_{n+1}}{q_{n+1}-r_{n+1}}
 =\frac{-\alpha_nr_n(3+\eta_n)}{-2\alpha_nr_n\eta_n}
 =\frac{\eta_n+3}{2\eta_n}=v_n^4.
\]
This is precisely $E_{v_n}$. In particular, the point
$(r_{n+1},0)$ becomes $(0,0)$, while the other two points become
$(1,0)$ and $(v_n^4,0)$. The zero abscissa in the Legendre model
is therefore compatible with the nonzero abscissa imposed in
the root chart: the coordinates have been translated and scaled.
The map in Proposition~\ref{prop:legendre} below fixes an explicit
normalization of the quotient by this kernel, including the
choice of its $y$-coordinate.

We also recall the differential attached to these coordinates.
On the smooth projective curve $E_v$, the expression
\[
 \xi_v=\frac{dx}{2y}
\]
extends from the affine chart to a nowhere-vanishing global
regular differential. It is invariant under translations and
spans the one-dimensional space of regular differentials on
$E_v$. 
%Although the displayed denominator vanishes at points
%of order two, the differential is regular there: writing
%$f(x)=x(x-1)(x-v^4)$, the identity $2y\,dy=f'(x)\,dx$ gives
%$dx/(2y)=dy/f'(x)$, and $f'(x)\ne0$ at a root of $f$.
%At the point at infinity the usual local parameter $-x/y$
%likewise shows that this differential is regular and nonzero.

\begin{lemma}\label{lem:admissible}
Every successor of an admissible parameter is admissible.
\end{lemma}

\begin{proof}
If $w^2=\displaystyle \frac{v^2+1}{2v}$, then
\begin{equation}\label{eq:w4}
 w^4=\frac{(v^2+1)^2}{4v^2},\qquad
 w^4-1=\frac{(v^2-1)^2}{4v^2}.
\end{equation}
Thus $w=0$ would imply $v^2=-1$, and $w^4=1$ would imply $v^2=1$.
Both contradict admissibility.
\end{proof}

\begin{proposition}\label{prop:legendre}
For an admissible recurrence step $v\mapsto w$, the map
\begin{equation}\label{eq:psi}
 \Psi_v(x,y)=\left(
 \frac{(x+v^2)^2}{4v^2x},
 \frac{y(x^2-v^4)}{8v^3x^2}\right)
\end{equation}
is a $2$-isogeny $E_v\longrightarrow E_w$ with kernel
$\{O,(0,0)\}$. If $(x',y')$ denote the coordinates on $E_w$
and $\xi_w=dx'/(2y')$ its global regular differential, the
pullback to $E_v$ satisfies
\begin{equation}\label{eq:differential}
 \Psi_v^*\xi_w=2v\xi_v.
\end{equation}
Successive maps of this form do not backtrack.
\end{proposition}

\begin{proof}
Write $\lambda=v^4$ and $\lambda'=w^4$ for the Legendre parameters
of the source and target, respectively. Their equations, in their respective coordinates,
are
\[
 E_v:\ y^2=x(x-1)(x-\lambda),\qquad
 E_w:\ (y')^2=x'(x'-1)(x'-\lambda').
\]
The image point has coordinates
\[
 x'=\frac{(x+v^2)^2}{4v^2x},\qquad
 y'=\frac{y(x^2-v^4)}{8v^3x^2}.
\]
Using $\lambda'=(v^2+1)^2/(4v^2)$, we compute
\begin{align*}
 x'(x'-1)(x'-\lambda') &=\frac{(x+v^2)^2}{4v^2x}\cdot
\frac{(x-v^2)^2}{4v^2x} \cdot \frac{(x-1)(x-v^4)}{4v^2x} \\
 &=\frac{(x^2-v^4)^2(x-1)(x-v^4)}{64v^6x^3}=\frac{y^2(x^2-v^4)^2}{64v^6x^4}=(y')^2.
\end{align*}
Thus the formulas define a nonconstant rational map into $E_w$,
which extends to a morphism between the smooth projective curves.

The function $x'$ has poles exactly at $O$ and $(0,0)$ on $E_v$,
of order two at each point. Indeed, $x$ has a pole of order two
at $O$ and a zero of order two at $(0,0)$, and the numerator
$(x+v^2)^2$ is nonzero at the latter point. Since the coordinate
function on $E_w$ has its only pole at the point at infinity,
both of these points map to that point. Comparing pole degrees
gives $2\deg\Psi_v=4$, so $\deg\Psi_v=2$.
A nonconstant morphism of elliptic curves sending the origin
to the origin is an isogeny; hence its kernel is exactly
$\{O,(0,0)\}$. The coordinate change~\eqref{eq:roottolegendre}
explains why the kernel generator has abscissa zero in this model.

For the differential, we compute the pullback of the form on
the target by substituting the coordinate functions $x'$ and $y'$.
Differentiation gives
\[
 \frac{dx'}{dx}=\frac{x^2-v^4}{4v^2x^2};
\]
therefore
\[
 \Psi_v^*\xi_w
 =\frac{d(x'(x,y))}{2y'(x,y)}
 =\frac{\dfrac{x^2-v^4}{4v^2x^2}\,dx}
 {2\,\dfrac{y(x^2-v^4)}{8v^3x^2}}
 =2v\frac{dx}{2y}=2v\xi_v.
\]
This calculation on an affine open set proves the identity of
global differentials in~\eqref{eq:differential}.

Finally, substitution shows that the remaining two nonzero
points of order two, $(1,0)$ and $(\lambda,0)$, both map to $(\lambda',0)$.
Their common image generates the kernel of the dual isogeny.
The next map $\Psi_w$ has kernel generated by $(0,0)$ on $E_w$;
these points differ because $\lambda'\ne0$. This proves that successive
maps do not backtrack.
\end{proof}

%%%%%%%%%%%%%%%%%%%%%%%%%%%%%
\section{The Arithmetic Geometric Mean}\label{sec:agm}
%%%%%%%%%%%%%%%%%%%%%%%%%%%%%
The recurrence for $v_n$ can also be viewed as a normalized form
of the arithmetic-geometric mean, or AGM. We develop this relation
in both directions, keeping track of the choices of square roots.
The classical AGM is discussed in~\cite{Cox}, and its $2$-adic
analogue in~\cite{REF-2}. Finite-field versions and their relation
to isogenies are studied in~\cite{REF-1,AGMaquariums}.

%%%%%%%%%%%%%%%%%%%%%%%%%%%%%
\subsection{The AGM recurrence}\label{sec:agmrecurrence}
%%%%%%%%%%%%%%%%%%%%%%%%%%%%%
Let $K$ be a field of characteristic different from $2$.
Starting with $a_0,b_0\in K$, an AGM sequence is obtained by
iterating
\begin{equation}\label{eq:agmdefinition}
 a_{n+1}=\frac{a_n+b_n}{2},\qquad
 b_{n+1}^2=a_nb_n.
\end{equation}
The first equation takes the arithmetic mean, whereas the
second chooses a geometric mean. Over a general field there
is no preferred sign for this square root. If $a_nb_n$ is a
nonzero square in $K$, the step has two possible successors,
\[
 (a_n,b_n)\longmapsto
 \left(\frac{a_n+b_n}{2},\sqrt{a_nb_n}\right),\qquad
 \left(\frac{a_n+b_n}{2},-\sqrt{a_nb_n}\right).
\]
If it is not a square in $K$, the branch cannot be continued
within $K$, although it can be continued over an algebraic
extension. Allowing both roots at each step therefore produces
a directed graph rather than a single sequence.

For the connection with elliptic curves, we use pairs satisfying
\[
 a_nb_n\ne0,\qquad a_n\ne\pm b_n.
\]
These conditions are preserved whenever an AGM successor exists.
Indeed, $a_{n+1}\ne0$ because $a_n\ne-b_n$, and $b_{n+1}\ne0$
because $a_nb_n\ne0$. Moreover,
\[
 a_{n+1}^2-b_{n+1}^2
 =\frac{(a_n+b_n)^2}{4}-a_nb_n
 =\frac{(a_n-b_n)^2}{4}\ne0.
\]
Thus a branch starting from such a pair stays away from the
singular parameters of the associated Legendre curves.

%%%%%%%%%%%%%%%%%%%%%%%%%%%%%
\subsection{From recurrence parameters to AGM pairs}\label{sec:agmforward}
%%%%%%%%%%%%%%%%%%%%%%%%%%%%%
\begin{proposition}\label{prop:agm}
For any admissible recurrence branch, set
\[
 b_0=1,\qquad b_n=\prod_{i=0}^{n-1}v_i,\qquad a_n=v_n^2b_n.
\]
Then
\begin{equation}\label{eq:agm}
 a_{n+1}=\frac{a_n+b_n}{2},\qquad b_{n+1}^2=a_nb_n.
\end{equation}
Conversely, an AGM sequence with nonzero entries yields
$v_n=b_{n+1}/b_n$ and $v_n^2=a_n/b_n$, satisfying the recurrence
wherever these terms are defined. If also $a_n\ne\pm b_n$,
the resulting parameters are admissible.
\end{proposition}
\begin{proof}
The empty product gives $b_0=1$ and hence $a_0=v_0^2$.
From the definition of $b_n$ we have $b_{n+1}=v_nb_n$, so
\[
 b_{n+1}^2=v_n^2b_n^2=(v_n^2b_n)b_n=a_nb_n.
\]
For the arithmetic mean, use the recurrence for $v_n$ to obtain
\begin{align*}
 a_{n+1}
 &=v_{n+1}^2b_{n+1}
 =\frac12\left(v_n+\frac1{v_n}\right)v_nb_n\\
 &=\frac{v_n^2+1}{2}b_n
 =\frac{v_n^2b_n+b_n}{2}
 =\frac{a_n+b_n}{2}.
\end{align*}
All entries are nonzero. Also $a_n/b_n=v_n^2$, and admissibility
gives $(a_n/b_n)^2=v_n^4\ne1$, so $a_n\ne\pm b_n$.

Conversely, suppose that $(a_n,b_n)$ is an AGM sequence with
nonzero entries, and define
\[
 v_n=\frac{b_{n+1}}{b_n}.
\]
The geometric-mean relation shows directly that
\[
 v_n^2=\frac{b_{n+1}^2}{b_n^2}
 =\frac{a_nb_n}{b_n^2}=\frac{a_n}{b_n}.
\]
Thus this choice of square root of $a_n/b_n$ is fixed by the
chosen next term $b_{n+1}$. Using both AGM relations, we obtain
\begin{align*}
 v_{n+1}^2
 &=\frac{a_{n+1}}{b_{n+1}}
 =\frac{a_n+b_n}{2b_{n+1}}\\
 &=\frac{v_n^2b_n+b_n}{2v_nb_n}
 =\frac{v_n^2+1}{2v_n}
 =\frac12\left(v_n+\frac1{v_n}\right).
\end{align*}
Finally, if $a_n\ne\pm b_n$, then
$v_n^4=(a_n/b_n)^2\ne1$, proving admissibility.
\end{proof}

The normalization $b_0=1$ is convenient but not essential.
For any $c\ne0$, multiplying every pair by $c$ preserves the
AGM equations, since
\[
 ca_{n+1}=\frac{ca_n+cb_n}{2},\qquad
 (cb_{n+1})^2=(ca_n)(cb_n).
\]
For a general nonzero initial value $b_0$, the inverse construction
therefore recovers
\begin{equation}\label{eq:agmproducts}
 b_n=b_0\prod_{i=0}^{n-1}v_i,\qquad a_n=v_n^2b_n.
\end{equation}
In particular, the product of the recurrence parameters records
the change in the geometric-mean coordinate:
\[
 \prod_{i=0}^{N-1}v_i=\frac{b_N}{b_0}.
\]
This also describes the sign choices without introducing new
square roots: choosing the other value of $b_{n+1}$ changes the
sign of $v_n=b_{n+1}/b_n$.

%%%%%%%%%%%%%%%%%%%%%%%%%%%%%
\subsection{Associated curves and the quadratic field}\label{sec:agmcurves}
%%%%%%%%%%%%%%%%%%%%%%%%%%%%%
To an AGM pair $(a,b)$ with $ab\ne0$ and $a\ne\pm b$, associate
the Legendre curve
\[
 \mathcal E_{(a,b)}:\quad
 y^2=x(x-1)(x-\lambda),\qquad \lambda=\left(\frac ab\right)^2.
\]
If $(a,b)=(cv^2,c)$ as above, then
\[
 \lambda=\left(\frac{cv^2}{c}\right)^2=v^4,\qquad
 \mathcal E_{(cv^2,c)}=E_v.
\]
Thus the curve attached to the AGM pair is exactly the model
of Section~\ref{sec:legendre}, with invariant $J(v)$.
If $(a',b')$ is an AGM successor, then its parameter is
\[
 \lambda'=\left(\frac{a'}{b'}\right)^2
 =\frac{(a+b)^2}{4ab}
 =\frac{(v^2+1)^2}{4v^2}=w^4,
 \qquad w^2=\frac12(v+v^{-1}).
\]
Consequently the $2$-isogeny $\Psi_v:E_v\to E_w$ of
Proposition~\ref{prop:legendre} directly relates the curves
associated to consecutive AGM pairs, with no further change
of coordinates.

Over $k=\FF_{p^2}$, Proposition~\ref{prop:agm} identifies the
field-of-definition question for these sequences with that for
the recurrence parameters. Products and quotients of nonzero
elements remain in $k$, so a branch of $v_n$ in $k$ gives AGM
pairs in $k$, and an AGM branch in $k$ gives its parameters
$v_n=b_{n+1}/b_n$ in $k$.

Theorem~\ref{thm:A}, proved in the next section, will in particular
give the sufficient condition from this viewpoint: if
$v_0\in k$ is admissible and $J(v_0)$ is supersingular, every
AGM branch starting at $(v_0^2,1)$ remains in $k$.
Indeed, the first geometric mean is either $v_0$ or $-v_0$.
Both parameters have the same supersingular invariant, and the
theorem applies to all subsequent choices from either one.
The proof of Theorem~\ref{thm:A} uses the algebraic correspondence
above and the independently cited structure of AGM graphs;
it does not assume this rationality conclusion in advance.

%%%%%%%%%%%%%%%%%%%%%%%%%%%%%
\section{Theorem A: rational branches and supersingularity}\label{sec:A}
%%%%%%%%%%%%%%%%%%%%%%%%%%%%%
We use the AGM graph on
\[
 V(k)=\{(a,b)\in(k^\times)^2:a\ne\pm b\},
\]
with arrows $(a,b)\mapsto((a+b)/2,\pm\sqrt{ab})$ when these roots
belong to $k$. As in Section~\ref{sec:agmcurves}, we associate
the Legendre parameter $\lambda=(a/b)^2$ to a vertex $(a,b)$.
The reference~\cite{AGMaquariums} uses the reciprocal parameter
$\lambda_{\mathrm{lit}}=(b/a)^2$. These conventions give curves
isomorphic over $k$: putting $t=a/b\in k^\times$, the map
\[
 (x,y)\longmapsto\left(\frac{x}{t^2},\frac{y}{t^3}\right)
\]
sends $y^2=x(x-1)(x-t^2)$ to
$y^2=x(x-1)(x-t^{-2})$.
The AGM graph itself is unchanged, and so is the ordinary or
supersingular nature of the curve at each vertex. The following
results therefore apply with our convention as well.
We recall two results from~\cite[Lemma 5.1 and Proposition 4.3]{AGMaquariums},
using the numbering of its preprint version:
\begin{enumerate}
\item A vertex associated with a supersingular curve has both
children in $V(\FF_{p^2})$.
\item A strongly connected component containing a directed cycle
and associated with ordinary curves is a simple directed cycle.
Here a \emph{strongly connected component} means a maximal set
of vertices such that any two vertices can be joined by directed
paths in both directions.
An elliptic curve $E$ in characteristic $p$ is \emph{ordinary}
if it has exactly $p$ geometric points killed by multiplication
by $p$, that is, $\#E[p](\overline{k})=p$; in contrast, a
supersingular elliptic curve has only the identity point in
$E[p](\overline{k})$.
A \emph{simple directed cycle} means a cyclic sequence of distinct
vertices whose arrows go from each vertex to the next and from
the last back to the first, with no other arrows inside the
component. Thus each vertex has exactly one incoming and one
outgoing arrow within the component; arrows entering or leaving
the component are not excluded.
\end{enumerate}
These results will be used explicitly, rather than assuming that
rationality of $j$ alone implies rationality of its root data.
For the first result, the underlying fact is that every supersingular
Legendre parameter is a fourth power in $\FF_{p^2}$;
see also~\cite[Proposition 3.1]{AuerTop2002}.

\begin{maintheorem}[Rational-branch characterization]\label{thm:A}
Let $p\geq5$ be prime, $k=\FF_{p^2}$, and let $v_0\in k$ be
admissible. Then $E_{v_0}$ is supersingular if and only if every
finite branch of~\eqref{eq:recurrence} starting from $v_0$ belongs
to $k$. Equivalently, every required square root, for every
successive choice of signs, lies in $k$.
\end{maintheorem}

\begin{proof}
Suppose first that $E_v$ is supersingular and $v\in k$ is
admissible. Put $r=(v+v^{-1})/2$. By Proposition~\ref{prop:legendre},
the curve with Legendre parameter $r^2$ is supersingular.
By~\eqref{eq:w4}, $r\ne0,\pm1$, so $(r,1)\in V(k)$,
whose associated Legendre parameter is $r^2$.
The first AGM result recalled above implies that $r$ has its two
square roots in $k$. Both successors of $v$ therefore lie in $k$.
They are admissible by Lemma~\ref{lem:admissible}, and their curves
are supersingular by isogeny invariance. Induction proves the assertion.

Conversely, suppose every branch from $v_0$ stays in $k$.
Let $S$ consist of all parameters occurring at positive depth in
these branches. This is a finite nonempty set, it is closed under
negation, meaning that $v\in S$ implies $-v\in S$, and it contains
both successors of each of its elements.
Consider the finite subset of AGM vertices
\[
 W=\{(cv^2,c):c\in k^\times,\ v\in S\}\subseteq V(k).
\]
It is closed under both AGM arrows. Indeed, for $(a,b)=(cv^2,c)$
the children have second coordinate $\pm cv$. For the plus sign,
put $c'=cv$ and choose a successor $w$ of $v$. Then
\[
 (a',b')=\left(c\frac{v^2+1}{2},cv\right)=(c'w^2,c').
\]
For the minus sign use $-v\in S$ instead. In both cases the child
belongs to $W$ and the two children are distinct.

We now choose a subset of $W$ as follows. For each vertex $P\in W$,
let $C(P)$ consist of $P$ and all vertices that can be reached
from $P$ by following finitely many arrows. All these vertices
belong to $W$, because both successors of every vertex of $W$
remain in $W$. Since $W$ is finite and nonempty, we can choose
$P_0\in W$ for which $C(P_0)$ has the smallest possible number
of vertices. Put $C=C(P_0)$.

If $Q\in C$, every vertex reachable from $Q$ is also reachable
from $P_0$, so $C(Q)\subseteq C$. By the minimal choice of $P_0$,
$C(Q)$ cannot have fewer vertices than $C$, and hence $C(Q)=C$.
Thus, starting at any vertex of $C$, we can reach every other
vertex of $C$ by following arrows. Moreover, both successors
of each vertex of $C$ lie in $C$: an arrow to either successor
can simply be appended to a path from $P_0$.

No vertex outside $C$ can be reached from $P_0$, even in the
full AGM graph, since every path starting in $W$ stays in $W$.
Consequently $C$ cannot be enlarged while preserving the property
that every vertex can reach every other one. This is precisely
the strongly connected component to which the second AGM result
applies. It contains a directed cycle: following successive arrows
inside the finite set $C$ must eventually repeat a vertex.
All associated elliptic curves are isogenous to
$E_{v_0}$: the curve attached to $(cv^2,c)$ is exactly $E_v$,
and $v$ is a descendant
of $v_0$. If $E_{v_0}$ were ordinary, the second AGM result would
force this component to be a simple directed cycle. A simple cycle
has only one outgoing edge within the component at each vertex,
contradicting the two distinct children. Hence $E_{v_0}$ is supersingular.
\end{proof}

\begin{example}\label{ex:seed}
	Let $v_0$ satisfy $v_0^4=2$. It is admissible, and
	$J(v_0)=256\,(2^2-2+1)^3/(2^2-2)^2=1728$, so $E_{v_0}\colon y^2=x(x-1)(x-2)$
	is supersingular exactly when $p\equiv 3 \pmod 4$. In that case $(p^2-1)/4$
	is a multiple of $p-1$, so every nonzero element of $\mathbb{F}_p$ is a
	fourth power in $k$; hence $v_0\in k$, and Theorem~\ref{thm:A} shows that every branch
	from $v_0$ remains in $k$. By contrast, if $p\equiv 1 \pmod 8$ then $2$ is
	still a fourth power in $k$, so $v_0\in k$, but $E_{v_0}$ is ordinary and
	Theorem~\ref{thm:A} shows that some branch from $v_0$ leaves $k$. In the notation of
	Proposition~\ref{prop:agm}, $v_0$ corresponds to the AGM seed
	$(a_0,b_0)=(\sqrt{2},1)$, which over $\mathbb{R}$ is, up to scaling, the
	starting point of the Gauss--Legendre algorithm for $\pi$; see~\cite{Cox}.
\end{example}

%%%%%%%%%%%%%%%%%%%%%%%%%%%%%
\subsection{Three Frobenius-stable directions}
%%%%%%%%%%%%%%%%%%%%%%%%%%%%%
The following complementary results express a finite rationality
criterion directly in the Tate module. They do not require a choice
of root chart or a convention for how many levels that chart records.

\begin{lemma}\label{lem:cylinder}
Let $T$ be free of rank two over $\ZZ_2$. If $A\in\End_{\ZZ_2}(T)$
stabilizes every rank-one direct summand in a nonempty open subset
of $\PP^1(\ZZ_2)$, then $A$ is scalar.
\end{lemma}
\begin{proof}
After a change of basis, let $e_1,e_2$ denote the resulting
$\ZZ_2$-basis of $T$, chosen so that the subset contains
$L_u=\langle e_1+2^ru e_2\rangle$, $u\in\ZZ_2$, for some $r\geq1$.
In this basis, write $A=\left(\begin{smallmatrix}a&b\\c&d\end{smallmatrix}\right)$.
Stability at $u=0$ gives $c=0$, and at $u=1,2$ gives
\[
 2^r(d-a)-2^{2r}b=0,\qquad
 2^{r+1}(d-a)-2^{2r+2}b=0.
\]
Subtracting twice the first equation yields $2^{2r+1}b=0$.
Since $\ZZ_2$ is torsion-free, $b=0$ and then $d=a$.
\end{proof}

\begin{corollary}\label{cor:scalar}
Let $E/k$ be an elliptic curve. If the $p^2$-power Frobenius
stabilizes a nonempty open subset of lines of $T_2(E)$, then its
characteristic polynomial is $(X-p)^2$ or $(X+p)^2$, and $E$ is
supersingular.
\end{corollary}
\begin{proof}
By Lemma~\ref{lem:cylinder} its action is $aI$. The determinant
is $p^2$, so $a^2=p^2$ in $\ZZ_2$ and $a=\pm p$.
The resulting trace $\pm2p$ implies supersingularity.
\end{proof}

\begin{theorem}[Finite three-direction criterion]\label{thm:finite}
Let $E/k$ be an elliptic curve and let $P,Q$ be a basis of
$E[2^n]$ over $\ZZ/2^n\ZZ$. Let $n\geq2r+2$, $r\geq1$.
Assume that the three subgroups
\begin{equation}\label{eq:directions}
 \langle P\rangle,\quad \langle P+2^rQ\rangle,\quad
 \langle P+2^{r+1}Q\rangle
\end{equation}
are stable under the $p^2$-power Frobenius. If
\begin{equation}\label{eq:depth}
 2^{2(n-r-1)}>4p^2,
\end{equation}
then $E$ is supersingular.
\end{theorem}
\begin{proof}
Lift the chosen basis to $T_2(E)$ and write the Frobenius matrix
as $A=\left(\begin{smallmatrix}a&b\\c&d\end{smallmatrix}\right)$
with entries in $\ZZ_2$. Stability gives
\begin{align*}
 c&\equiv0\pmod{2^n},\\
 d-a-2^rb&\equiv0\pmod{2^{n-r}},\\
 d-a-2^{r+1}b&\equiv0\pmod{2^{n-r-1}}.
\end{align*}
Subtracting yields $2^{n-2r-1}\mid b$ and
$2^{n-r-1}\mid(d-a)$. Therefore, for its integer trace $t$,
\[
 \Delta=t^2-4p^2=(a-d)^2+4bc
 \quad\hbox{satisfies}\quad 2^{2(n-r-1)}\mid\Delta.
\]
Here $4bc$ is divisible by $2^{2n-2r+1}$, so the precision of the
matrix reduction suffices for the asserted divisibility.
Hasse's bound gives $|\Delta|\leq4p^2$; hence~\eqref{eq:depth}
forces $\Delta=0$, and $E$ is supersingular.
\end{proof}

For $r=1$, the inequality $n>\log_2p+3$ suffices.
Thus three finite chains of length $n$, rather than infinite
chains, suffice to prove supersingularity, provided that they
start at the same curve $E$, consist of $2$-isogenies defined
over $k=\FF_{p^2}$, do not backtrack, and that the kernels of
their composite isogenies are
$\langle P\rangle$, $\langle P+2Q\rangle$, and
$\langle P+4Q\rangle$ for a common basis $P,Q$ of $E[2^n]$.
Indeed, these kernels are cyclic and Frobenius-stable, so the
theorem applies. However, the criterion does not provide a way
to choose three chains with these kernels among the branches
of the recurrence graph generated by a given parameter.
To use it in that graph, one must still identify suitable chains
and verify the required relation between their kernels.
The existence of three arbitrary distinct chains of the stated
length is therefore not sufficient for this application.

%%%%%%%%%%%%%%%%%%%%%%%%%%%%%
\section{The edge parameter and exceptional configurations}\label{sec:graph}
%%%%%%%%%%%%%%%%%%%%%%%%%%%%%
We consider geometric isomorphism classes and count outgoing
isogenies by their kernels. Each curve has three order-two
subgroups. At exceptional automorphism loci, identifying maps or
identifying dual arrows can change the graphical multiplicity;
the kernel convention remains the one used in proofs.
For background on loops and multiple edges see~\cite{Ghantous,REF-5}.

The classical modular polynomial of level $2$ is the symmetric
polynomial $\Phi_2(X,Y)\in\ZZ[X,Y]$ whose vanishing characterizes
pairs of $j$-invariants of elliptic curves related by a $2$-isogeny.
Here $X$ is the invariant of the source and $Y$ that of the target.
More precisely, after reduction to our characteristic $p\geq5$,
\[
 \Phi_2(j(E),Y)=\prod_{H\subset E[2],\ |H|=2}
 \bigl(Y-j(E/H)\bigr).
\]
The product runs over the three order-two subgroups over $\bar k$.
It is monic of degree three in $Y$, and repeated roots record
distinct kernels whose quotients have the same invariant.
This describes geometric isogenies; it does not by itself assert
that an isogeny between prescribed models is defined over $k$.

The standard coefficient formula, tabulated
in~\cite{SutherlandModularPolynomials}, is
\begin{align*}
 \Phi_2(X,Y)={}&X^3+Y^3-X^2Y^2+1488XY(X+Y)\\
 &-162000(X^2+Y^2)+40773375XY\\
 &+8748000000(X+Y)-157464000000000.
\end{align*}
Using $X^3+Y^3=(X+Y)^3-3XY(X+Y)$ and
$X^2+Y^2=(X+Y)^2-2XY$, we obtain the equivalent form
\begin{align}\label{eq:Phi2}
 \Phi_2(X,Y)={}&(X+Y)^3-X^2Y^2+1485XY(X+Y)\nonumber\\
 &-162000(X+Y)^2+41097375XY\nonumber\\
 &+8748000000(X+Y)-157464000000000.
\end{align}
To explain the edge parameter, translate the chosen kernel point
to $(0,0)$ and write the source as
\[
 y^2=x(x^2+Ax+B),\qquad B(A^2-4B)\ne0.
\]
Its quotient by $\langle(0,0)\rangle$ has a model
$y^2=x(x^2-2Ax+A^2-4B)$, obtained from the quotient formulas
of Section~\ref{sec:prelim} by a translation on the target.
Computing the two invariants and putting $t=16(A^2-4B)/B\ne0$
gives
\[
 X=\frac{256(A^2-3B)^3}{B^2(A^2-4B)}
   =\frac{(t+16)^3}{t},\qquad
 Y=\frac{16(A^2+12B)^3}{B(A^2-4B)^2}
   =\frac{(t+256)^3}{t^2}.
\]
Thus $t$ describes the isogeny with its chosen kernel.
Define its \emph{edge parameter} by
\[
 g=\frac{(t+16)(t+256)}{t}=t+272+\frac{4096}{t}.
\]
Replacing $t$ by $4096/t$ interchanges $X$ and $Y$ and leaves
$g$ unchanged, so this label does not distinguish an isogeny
from its dual. Direct multiplication and addition now give
\begin{equation}\label{eq:edge}
 XY=g^3,\qquad X+Y=g^2-495g+54000.
\end{equation}
In particular, $g$ is not an arbitrary choice of a cube root of
$XY$: it satisfies both relations in~\eqref{eq:edge}.
Substitution verifies $\Phi_2(X,Y)=0$, and elimination of $g$
from these two relations recovers the modular equation.
We use $g$ as a label for the unordered pair of endpoints of
an isogeny, retaining the possible coincidences discussed below.
Thus the two endpoints are roots of
\begin{equation}\label{eq:endpoints}
 Z^2-(g^2-495g+54000)Z+g^3=0.
\end{equation}
For a fixed endpoint $j$, setting $Z=j$ in this equation gives
the cubic whose roots are the incident parameter values,
counted with multiplicity:
\begin{equation}\label{eq:edgecubic}
 g^3-jg^2+495jg+j^2-54000j=0.
\end{equation}
Their sum is $j$. Distinct kernels need not give distinct values
of this parameter at a ramification point.

%%%%%%%%%%%%%%%%%%%%%%%%%%%%%
\subsection{Loops}
%%%%%%%%%%%%%%%%%%%%%%%%%%%%%
Substituting $X=Y=j$ into the modular equation yields
\[
 \Phi_2(j,j)=-(j-8000)(j-1728)(j+3375)^2.
\]
The discriminant of~\eqref{eq:edgecubic} as a polynomial in $g$ is
\[
 4j^2(j-1728)(j+3375)^2.
\]
For $j=8000$ the edge cubic factors as
\[
 (g-400)(g^2-7600g+920000),
\]
and the value $g=400$ gives $(Z-8000)^2$ in~\eqref{eq:endpoints}.
The diagrams record the corresponding algebraic configurations.
They must be read after reduction in the chosen characteristic;
for example, $8000=1728$ in characteristic seven.
\begin{figure}[htbp]
\centering\resizebox{.7524\linewidth}{!}{%
% Exact edge labels checked by revision/auditoria_trazas_figuras.py.
\begin{tikzpicture}[every node/.style={font=\scriptsize},
 vertex/.style={draw,circle,minimum size=13mm,inner sep=2pt,font=\small},
 every path/.style={draw=black,line width=.6pt}]
\node[vertex] (a) at (0,0) {$1728$};
\node[vertex] (b) at (3.3,0) {$287496$};
\path (a) edge[loop above,looseness=6] node[above] {$144$} (a);
\path (a) edge node[above] {$792$} (b);
\draw (b) -- (5,.7) node[right] {$143352+101574\sqrt2$};
\draw (b) -- (5,-.7) node[right] {$143352-101574\sqrt2$};
\node[vertex] (c) at (0,-3) {$-3375$};
\node[vertex] (d) at (3.3,-3) {$16581375$};
\path (c) edge[loop above,looseness=6] node[above] {$225$} (c);
\path (c) edge node[above] {$-3825$} (d);
\draw (d) -- (5,-2.3) node[right] {$8292600+3132675\sqrt7$};
\draw (d) -- (5,-3.7) node[right] {$8292600-3132675\sqrt7$};
\end{tikzpicture}
}
\caption{Distinct edge-parameter values at $j=1728$ and
$j=-3375$ and at their neighbours. Each label is a value of $g$;
the drawing does not encode kernel multiplicities.}\label{fig:loops}
\end{figure}
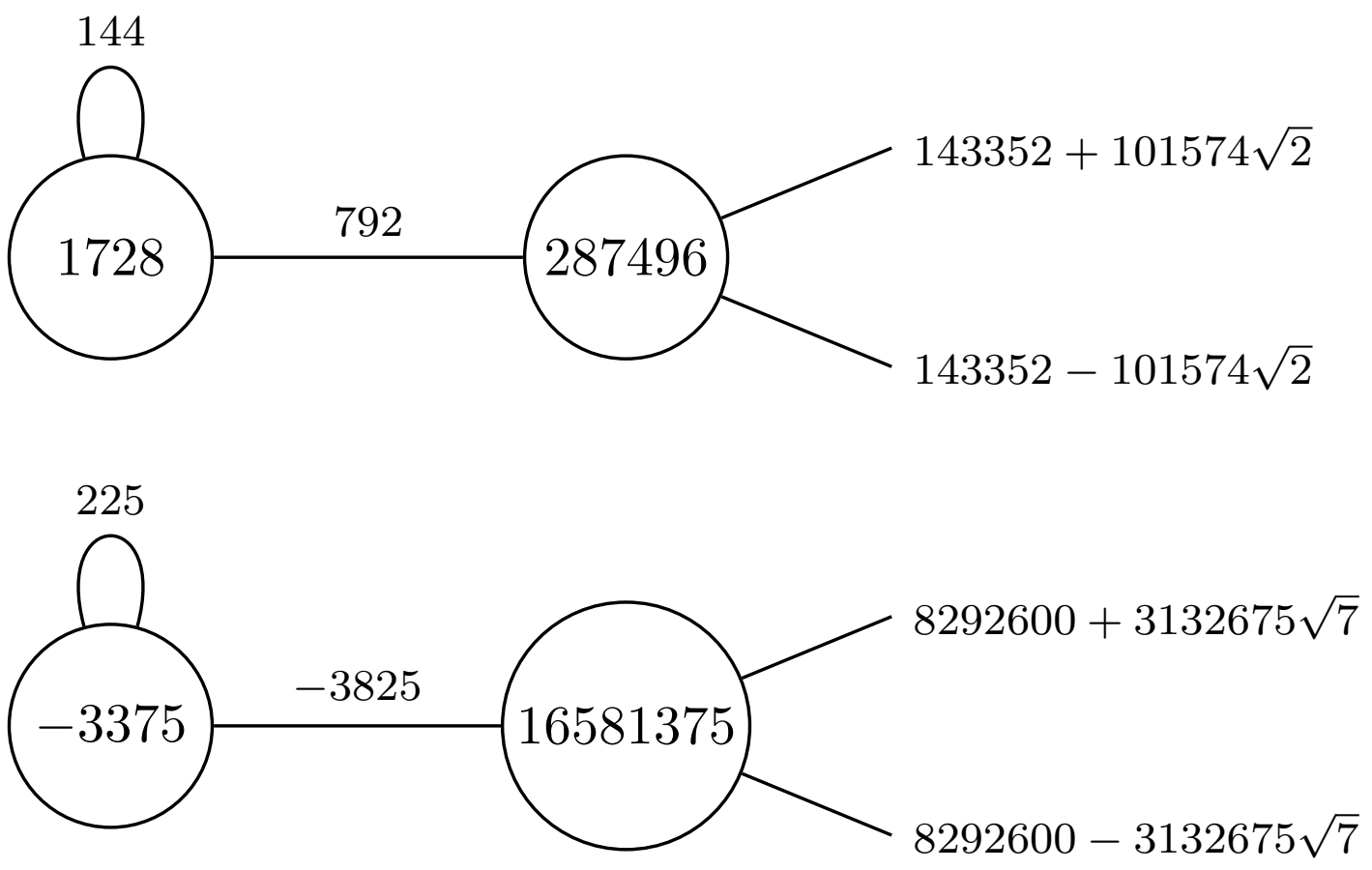
\begin{figure}[htbp]
\centering\resizebox{.95\linewidth}{!}{%
% Same checked labels as the original figure, in vector black.
\begin{tikzpicture}[every node/.style={font=\large},
 vertex/.style={draw,ellipse,minimum height=30mm,inner sep=5pt},
 every path/.style={draw=black,line width=.6pt}]
\node[vertex] (a) at (0,0) {$26125000+18473000\sqrt2$};
\node[draw,circle,minimum size=13mm,inner sep=2pt] (b) at (6,0) {$8000$};
\node[vertex] (c) at (12,0) {$26125000-18473000\sqrt2$};
\path (b) edge[loop above,looseness=6] node[above] {$400$} (b);
\path (a) edge node[above=48pt,xshift=-18pt] {$3800+2600\sqrt2$} (b);
\path (b) edge node[above=48pt,xshift=18pt] {$3800-2600\sqrt2$} (c);
\draw (a.north west) -- ++(-.7,.5);
\draw (a.south west) -- ++(-.7,-.5);
\draw (c.north east) -- ++(.7,.5);
\draw (c.south east) -- ++(.7,-.5);
\end{tikzpicture}
}
\caption{The configuration at $j=8000$, including the loop
parameter $g=400$.}\label{fig:8000}
\end{figure}
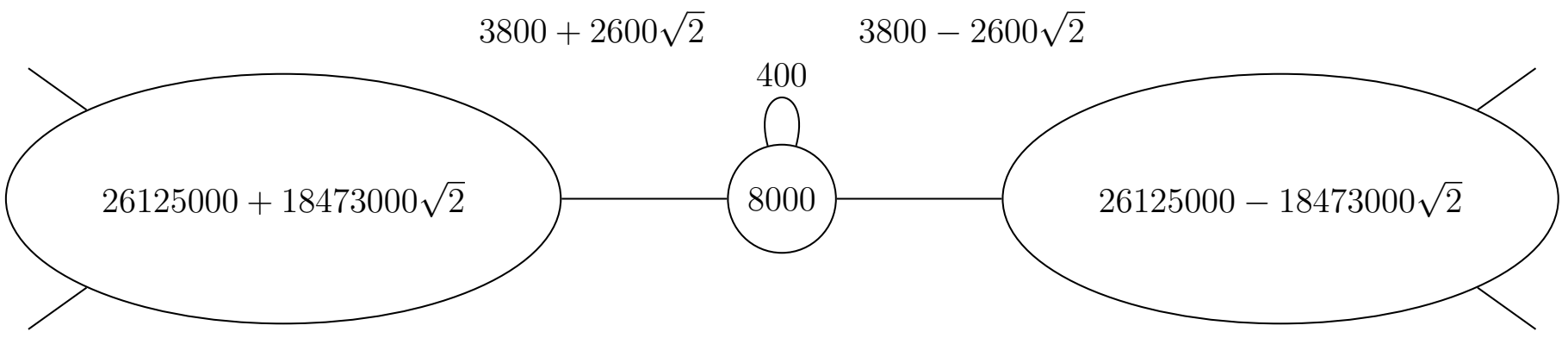

%%%%%%%%%%%%%%%%%%%%%%%%%%%%%
\subsection{Two distinct edge parameters with the same endpoints}
%%%%%%%%%%%%%%%%%%%%%%%%%%%%%
For endpoints $j_1,j_2$, the candidates are
\[
 g=\frac{495\pm\sqrt{29025+4(j_1+j_2)}}2.
\]
If two distinct candidates have the same cube $j_1j_2$, the
sum and product relations imply
\[
 j_1+j_2=-191025,\qquad j_1j_2=-495^3.
\]
Thus the endpoints satisfy
\begin{equation}\label{eq:double}
 Z^2+191025Z-121287375=0,
\end{equation}
and the parameters are $495(1\pm\sqrt{-3})/2$.
The requirement that these values and the relevant endpoints
remain distinct is essential after reduction.
Indeed the discriminants for the parameters and the endpoints are
$-3\cdot495^2$ and $5\cdot85995^2$, respectively; among primes
$p\geq5$, the displayed configuration degenerates at
$p\in\{5,7,11,13\}$.
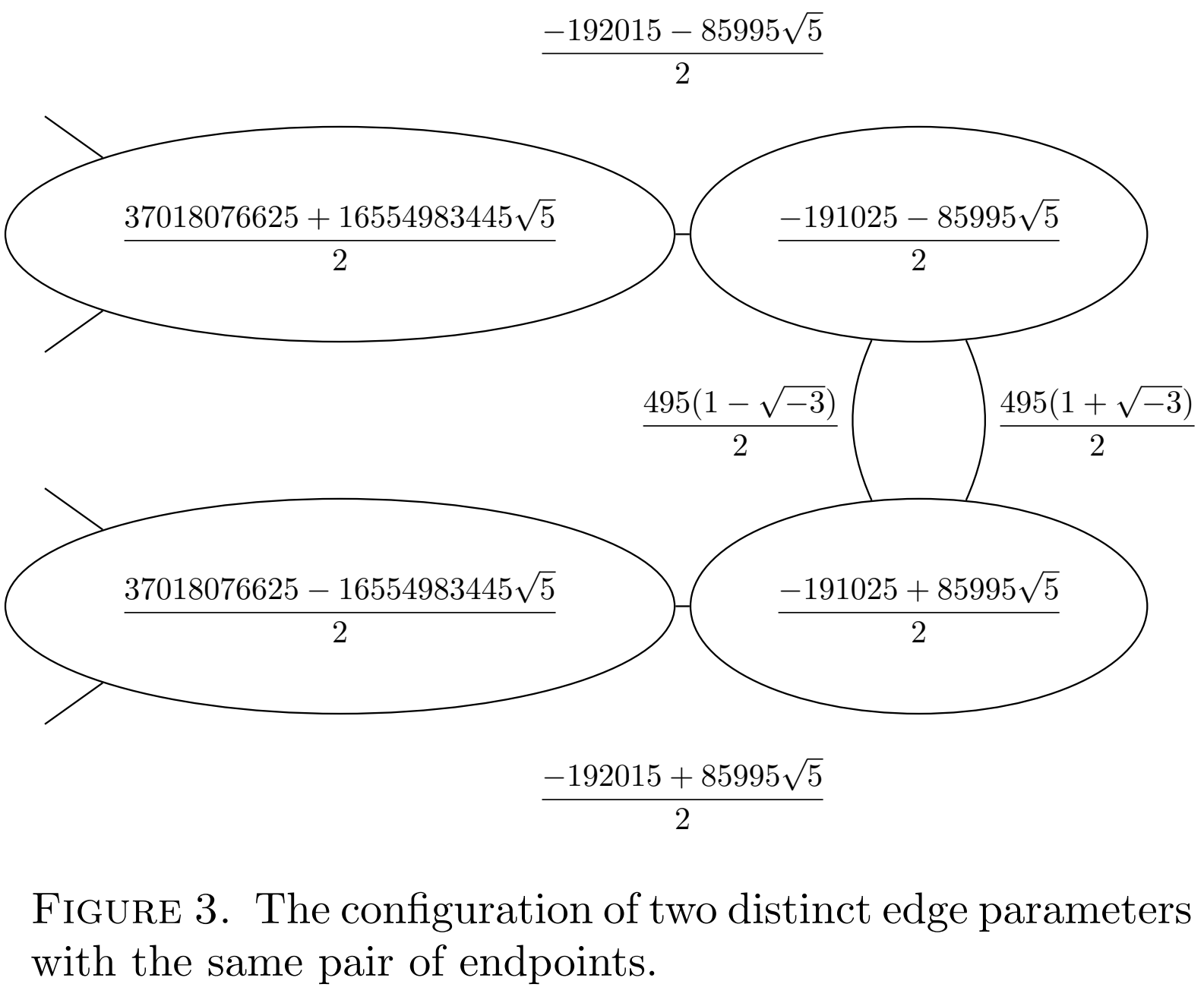
\begin{figure}[htbp]
\centering\resizebox{.836\linewidth}{!}{%
% Same checked labels and fractions as the original figure, in vector black.
\begin{tikzpicture}[every node/.style={font=\normalsize},
 vertex/.style={draw,ellipse,minimum height=26mm,inner sep=6pt},
 every path/.style={draw=black,line width=.6pt}]
\node[vertex] (a) at (0,0) {$\displaystyle\frac{37018076625+16554983445\sqrt5}{2}$};
\node[vertex] (b) at (7,0) {$\displaystyle\frac{-191025-85995\sqrt5}{2}$};
\node[vertex] (c) at (7,-4.5) {$\displaystyle\frac{-191025+85995\sqrt5}{2}$};
\node[vertex] (d) at (0,-4.5) {$\displaystyle\frac{37018076625-16554983445\sqrt5}{2}$};
\path (a) edge node[above=48pt] {$\displaystyle\frac{-192015-85995\sqrt5}{2}$} (b);
\path (b) edge[bend left=24] node[right] {$\displaystyle\frac{495(1+\sqrt{-3})}{2}$} (c);
\path (b) edge[bend right=24] node[left] {$\displaystyle\frac{495(1-\sqrt{-3})}{2}$} (c);
\path (c) edge node[below=48pt] {$\displaystyle\frac{-192015+85995\sqrt5}{2}$} (d);
\draw (a.north west) -- ++(-.7,.5);
\draw (a.south west) -- ++(-.7,-.5);
\draw (d.north west) -- ++(-.7,.5);
\draw (d.south west) -- ++(-.7,-.5);
\end{tikzpicture}
}
\caption{The configuration of two distinct edge parameters
with the same pair of endpoints.}\label{fig:double}
\end{figure}
\begin{figure}[htbp]
\centering\resizebox{.95\linewidth}{!}{%
% One label per unordered pair and g value; no kernel multiplicity implied.
\begin{tikzpicture}[every node/.style={font=\small},
 vertex/.style={draw,circle,minimum size=13mm,inner sep=2pt},
 every path/.style={draw=black,line width=.6pt}]
\node[vertex] (a) at (0,2.1) {$80$};
\node[vertex] (b) at (2.3,2.1) {$24$};
\node[vertex] (c) at (1.15,0) {$23$};
\node[vertex] (d) at (3.6,0) {$69$};
\node[vertex] (e) at (6,0) {$34$};
\node[vertex] (f) at (8.5,1.2) {$61+36\ima$};
\node[vertex] (g) at (8.5,-1.2) {$61+67\ima$};
\node[vertex] (h) at (11.6,1.2) {$20+94\ima$};
\node[vertex] (j) at (11.6,-1.2) {$20+9\ima$};
\path (a) edge[loop above,looseness=6] node[above] {$41$} (a);
\path (b) edge[loop above,looseness=6] node[above] {$19$} (b);
\path (d) edge[loop below,looseness=6] node[below] {$91$} (d);
\path (a) edge node[left] {$71$} (c);
\path (b) edge node[right] {$89$} (c);
\path (c) edge node[above] {$69$} (d);
\path (d) edge node[above] {$12$} (e);
\path (e) edge node[above,sloped] {$11+45\ima$} (f);
\path (e) edge node[below,sloped] {$11+58\ima$} (g);
\path (f) edge node[left] {$10$} (g);
\path (f) edge node[above] {$40+94\ima$} (h);
\path (g) edge node[below] {$40+9\ima$} (j);
\path (h) edge[bend right=22] node[left] {$90$} (j);
\path (h) edge[bend left=22] node[right] {$96$} (j);
\end{tikzpicture}
}
\caption{All supersingular vertices and distinct edge-parameter
values for $p=103$, with $\ima^2=-1$. Each unordered pair and each
value of $g$ is drawn once; kernel multiplicities are not encoded.}
\label{fig:p103}
\end{figure}
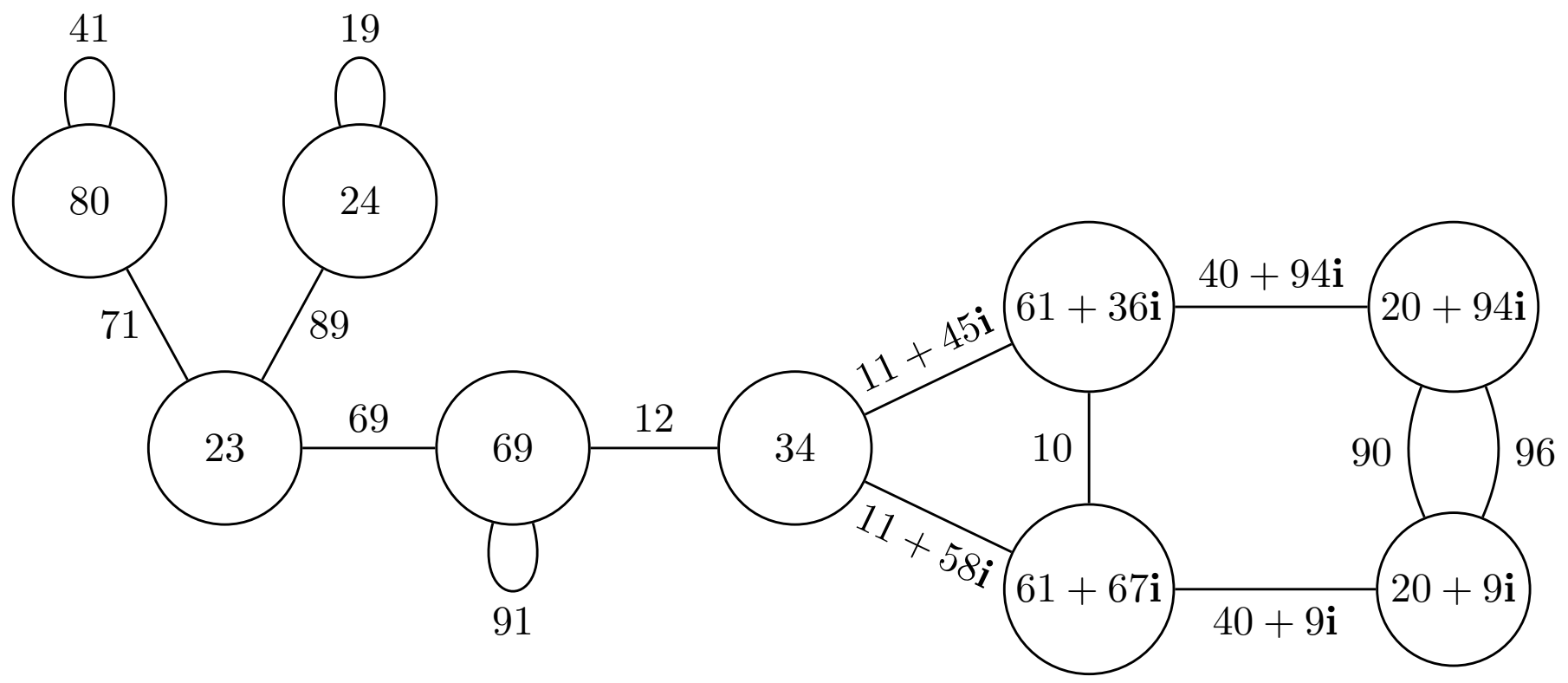

The exceptional $j$-values singled out by the preceding calculations are
\[
 0,\quad 1728,\quad -3375,\quad 8000,\quad
 \frac{-191025+85995\sqrt5}{2},\quad
 \frac{-191025-85995\sqrt5}{2}.
\]
More precisely, $j=0,1728,-3375$ are the values where the edge
cubic~\eqref{eq:edgecubic} has repeated roots; $j=1728,-3375,8000$
are the loop values, with respective loop parameters $g=144,225,400$;
and the last two displayed values are the endpoints associated with
the two distinct parameters $g=495(1\pm\sqrt{-3})/2$.
All values are interpreted after reduction in characteristic $p$,
so some may coincide; the last configuration degenerates at
$p=5,7,11,13$. In the supersingular graph, only those reductions
that are supersingular are retained. In particular, $j=0$ does
not occur there under the assumption $p\equiv7\pmod{12}$.

These configurations describe exceptional identifications in the
graph projected to $j$-invariants; they impose no additional
exclusions on admissible recurrence cycles in the parameters $v$.
In particular, Theorem~\ref{thm:B} remains valid for cycles passing
through these invariants, and does not require supersingularity.
The purpose of this section is to clarify the projected graph
and its multiplicities, rather than to supply a hypothesis for
that theorem. This does not mean that the parameter graph has
no cycles or coincidences of its own.

For completeness, the supersingularity conditions for the above
reductions, for primes $p\geq5$, are as follows:
$j=0$ requires $p\equiv2\pmod3$; $j=1728$ requires
$p\equiv3\pmod4$; $j=-3375$ requires $p=7$ or
$p\equiv3,5,6\pmod7$; $j=8000$ requires
$p\equiv5,7\pmod8$; and the two values involving $\sqrt5$
require $p=5$ or $p\equiv7,11,13,14\pmod{15}$.
These are the complex-multiplication reduction criteria for
discriminants $-3,-4,-7,-8,-15$, respectively; outside these
conditions the corresponding reductions are ordinary.
Thus, when $p\equiv7\pmod{12}$, $1728$ always occurs in the
supersingular graph and $0$ never does. Coincidences after
reduction must be distinguished from failure of supersingularity:
for example, at $p=7$, $8000=1728$ and this vertex is supersingular.

%%%%%%%%%%%%%%%%%%%%%%%%%%%%%
\section{Theorem B: cycle products and endomorphism traces}\label{sec:B}
%%%%%%%%%%%%%%%%%%%%%%%%%%%%%
The distinction between cycles of kernels, cycles of invariants
and their lifts is relevant to exact counting; see~\cite{Arpin2024,cryptoeprint:2025/155}.
The product identity below does not require a uniqueness theorem
for lifting arbitrary $j$-cycles.

%%%%%%%%%%%%%%%%%%%%%%%%%%%%%
\subsection{An order-four point fixing the sign}
%%%%%%%%%%%%%%%%%%%%%%%%%%%%%
\begin{definition}
An \emph{admissible recurrence cycle} of length $N$ is a tuple
$(v_0,\ldots,v_{N-1})$ of admissible parameters, meaning that
$v_n\ne0$ and $v_n^4\ne1$, satisfying
\[
 v_{n+1}^2=\frac12\left(v_n+\frac1{v_n}\right)
 \qquad (n=0,\ldots,N-1),\qquad v_N=v_0.
\]
Thus the last step returns to the initial parameter, and all
indices are read modulo $N$. The cycle is \emph{primitive}
if its least period is $N$.
The associated isogenies are the specific maps $\Psi_{v_n}$
in~\eqref{eq:psi}, not merely a list of $j$-invariants.
\end{definition}
For a recurrence cycle write
\begin{equation}\label{eq:ur}
 u_n=v_n^2,\qquad m_n=u_n^2,\qquad
 r_n=\frac{v_{n-1}^2-1}{2v_{n-1}}.
\end{equation}
Then $r_n^2=u_n^2-1$.
\begin{lemma}\label{lem:orderfour}
On $E_{v_n}$ the point
\begin{equation}\label{eq:R}
 R_n=(u_n^2+u_nr_n,\ u_nr_n(u_n+r_n))
\end{equation}
has exact order four, satisfies $2R_n=(m_n,0)$, and obeys
\begin{equation}\label{eq:transport}
 \Psi_{v_n}(R_n)=R_{n+1}.
\end{equation}
\end{lemma}
\begin{proof}
Write $u=u_n$, $r=r_n$, $s=u+r$. Then $s(u-r)=1$,
and for the point in~\eqref{eq:R} we have
\[
 x=us,\quad x-1=rs,\quad x-u^2=ur,\quad y=urs.
\]
This proves the curve equation. In the duplication formula for
$y^2=x^3-(1+u^2)x^2+u^2x$, the tangent slope at this point
is $s$. Indeed the numerator of the slope is $2ur s^2$ and
the denominator is $2ur s$. Hence
\[
 x(2R_n)=s^2+1+u^2-2us=u^2,
\]
and its $y$-coordinate is zero. Since $u,r,s$ are nonzero,
$R_n$ is not a point of order two.

For transport, put $v=v_n$ and use
\[
 u' =\frac{u+1}{2v},\quad r'=\frac{u-1}{2v},\quad u'+r'=v.
\]
The first coordinate of $\Psi_v(R_n)$ is $(u+1)/2=u'(u'+r')$.
Its second coordinate is
\[
 \frac{urs\,u^2(s^2-1)}{8v^3u^2s^2}
 =\frac{ur^2}{4v^3}=\frac{u^2-1}{4v}
 =u'r'(u'+r'),
\]
where $s^2-1=2rs$ and $u=v^2$ were used.
This gives~\eqref{eq:transport} with its stated sign.
\end{proof}

Recall that an endomorphism $\varphi$ of an elliptic curve $E$
acts linearly on the rank-two Tate module $T_\ell(E)$ for any
prime $\ell\ne p$. Its \emph{trace} $t=\operatorname{tr}(\varphi)$
is the trace of this linear action, that is, the sum of the two
diagonal entries of its matrix in any basis. It is an integer
independent of $\ell$ and of the chosen basis. The determinant
is $d=\deg\varphi$, so the characteristic polynomial and the
corresponding endomorphism identity are
\[
 X^2-tX+d,\qquad \varphi^2-[t]\varphi+[d]=0,
\]
where $[m]$ denotes multiplication by the integer $m$ on $E$.
Equivalently, the trace can be recovered from degrees by
$t=1+\deg\varphi-\deg(\varphi-\mathrm{id}_E)$, with the zero
endomorphism assigned degree zero. We apply this notion below
to the composite endomorphism obtained by going around the cycle.

\begin{maintheorem}[Product--trace identity]\label{thm:B}
Let $(v_0,\ldots,v_{N-1})$ be an admissible recurrence cycle
over $\bar k$, with $N\geq2$. Define
\[
 \Theta=\Psi_{v_{N-1}}\circ\cdots\circ\Psi_{v_0}\in\End(E_{v_0}),
 \qquad \omega=\prod_{n=0}^{N-1}v_n.
\]
Then $\deg\Theta=2^N$, and its integer trace $t$ satisfies
\begin{equation}\label{eq:tracecongruence}
 t\equiv1\pmod4,\qquad |t|\leq2^{1+N/2}.
\end{equation}
Moreover $\Theta(R_0)=R_0$ and
\begin{equation}\label{eq:B}
 \Theta^*\xi_{v_0}=2^N\omega\xi_{v_0},\qquad
 \boxed{2^N\omega^2-t\omega+1=0.}
\end{equation}
The last equation has coefficients in $\FF_p$, with
$\bar t=t\bmod p$, and its roots lie in $\FF_{p^2}$.
\end{maintheorem}
\begin{proof}
Degrees multiply. Lemma~\ref{lem:orderfour}, iterated around
the cycle, gives $\Theta(R_0)=R_0$. The integer trace and degree
satisfy the endomorphism equation
\[
 \Theta^2-[t]\Theta+[2^N]=0.
\]
Applying this equation to $R_0$ gives $[1-t+2^N]R_0=0$.
Since $R_0$ has exact order four and $N\geq2$, this proves
$t\equiv1\pmod4$. The usual positivity of the degree form gives
$t^2\leq4\deg\Theta$; see~\cite{Silverman}.

Multiplying~\eqref{eq:differential} around the cycle gives the
differential multiplier $\lambda=2^N\omega$. The differential
of the characteristic equation gives
\[
 \lambda^2-t\lambda+2^N=0.
\]
Substitute $\lambda=2^N\omega$ and divide by $2^N$, which is
nonzero in characteristic $p$, to obtain~\eqref{eq:B}.
\end{proof}

%%%%%%%%%%%%%%%%%%%%%%%%%%%%%
\subsection{The product polynomial}
%%%%%%%%%%%%%%%%%%%%%%%%%%%%%
\begin{corollary}\label{cor:productpolynomial}
Let
\[
 k_M=\max\{k\geq1:2k-1\leq2^{1+N/2}\}.
\]
Every admissible recurrence cycle of length $N\geq2$ satisfies
\begin{equation}\label{eq:productpolynomial}
 \prod_{k=1}^{k_M}
 \bigl(2^N\omega^2+(-1)^k(2k-1)\omega+1\bigr)=0.
\end{equation}
\end{corollary}
\begin{proof}
The integers congruent to one modulo four and satisfying the
trace bound are exactly $(-1)^{k+1}(2k-1)$ for $1\leq k\leq k_M$.
One of the factors therefore vanishes by Theorem~\ref{thm:B}.
\end{proof}

The polynomial~\eqref{eq:productpolynomial} has degree $2k_M$ in
$\omega$, but this does not imply that there are $2k_M$ primitive
cycles: different cycles can have the same product, not every
root need be realized on a specified curve, and a solution of
the period-$N$ equations may have smaller exact period.

\begin{remark}\label{rem:traceinteger}
	The product determines the trace modulo $p$, namely
	$t \bmod p = 2^N\omega+\omega^{-1}$. Since two integers congruent to
	$1$ modulo $4$ with the same residue modulo $p$ differ by a multiple of
	$4p$, the bound $|t|\le 2^{1+N/2}$ recovers the integer $t$ uniquely
	whenever $2^{N/2}<p$.
\end{remark}

%%%%%%%%%%%%%%%%%%%%%%%%%%%%%
\subsection{Conjugation and reversal}
%%%%%%%%%%%%%%%%%%%%%%%%%%%%%
Frobenius conjugation of every parameter sends a cycle over $k$
to another cycle, with product $\omega^p$.
For reversal put
\[
 T(v)=\frac{v+1}{v-1}.
\]
If $w^2=(v+v^{-1})/2$, then
\[
 \frac12\bigl(T(w)+T(w)^{-1}\bigr)
 =\frac{w^2+1}{w^2-1}=T(v)^2.
\]
Moreover $T(v)$ is admissible whenever $v$ is, since
\[
 T(v)^4-1=\frac{8v(v^2+1)}{(v-1)^4}\ne0,
\]
and $T(v)$ is defined and nonzero. Thus reversing the order and
applying $T$ gives an admissible recurrence cycle. Its product satisfies
\begin{equation}\label{eq:reciprocity}
 \omega_{T(\mathbf v)^\dagger}=\frac{2^{-N}}{\omega_{\mathbf v}}.
\end{equation}
Indeed,
\begin{align*}
 \prod_n T(v_n)
 &=\prod_n\frac{v_n^2-1}{(v_n-1)^2}
 =\prod_n\frac{(v_{n-1}-1)^2}
 {2v_{n-1}(v_n-1)^2}
 =\frac{2^{-N}}{\prod_n v_n}.
\end{align*}
This is the second root of the same product--trace quadratic.
In particular, reversal does not generally leave the product unchanged.

\begin{proposition}[The conjugate reverse cycle]\label{prop:reverseroot}
Let $\mathbf v=(v_0,\ldots,v_{N-1})$ be an admissible recurrence
cycle of length $N\geq2$, with product $\omega$ and integer trace
$t$ as in Theorem~\ref{thm:B}. Define its \emph{conjugate reverse
cycle}, with respect to $T(v)=(v+1)/(v-1)$, by
\[
 \mathbf v^{\mathrm{rev},T}
 =\bigl(T(v_{N-1}),T(v_{N-2}),\ldots,T(v_0)\bigr).
\]
This is an admissible recurrence cycle, and its product
$\omega^{\mathrm{rev},T}$ is the other root, counted with
multiplicity, of the product--trace equation
\[
 2^NX^2-tX+1=0.
\]
More precisely,
\[
 \omega^{\mathrm{rev},T}=\frac{1}{2^N\omega},\qquad
 2^NX^2-tX+1
 =2^N(X-\omega)(X-\omega^{\mathrm{rev},T}).
\]
These identities hold in characteristic $p$, with $t$ reduced
modulo $p$. The two roots coincide exactly when $t^2-2^{N+2}=0$
in that characteristic. Here conjugation refers to the map $T$,
not to Frobenius conjugation.
\end{proposition}
\begin{proof}
The preceding calculations show admissibility and give
$\omega^{\mathrm{rev},T}=2^{-N}/\omega$ by~\eqref{eq:reciprocity}.
Theorem~\ref{thm:B}, divided by the nonzero product $\omega$,
gives $t=2^N\omega+\omega^{-1}$.
Consequently the two products have sum $t/2^N$ and product
$2^{-N}$, proving the displayed factorization. The assertion
about coincident roots follows from its discriminant.
\end{proof}

\begingroup
\begin{remark}\label{rem:dualcycle}
The conjugate reverse cycle can be identified with the cycle of
dual isogenies, provided that the coordinate changes and the
initial vertex are chosen consistently. Write $\lambda_n=v_n^4$
and use $r_n$ from~\eqref{eq:ur}, with all indices modulo $N$.
The identity $T(v_n)^2=(v_{n+1}^2+1)/(v_{n+1}^2-1)$ and
\eqref{eq:w4} give
\[
 T(v_n)^4=\frac{\lambda_{n+2}}{\lambda_{n+2}-1}.
\]
Fix $i\in\bar k$ with $i^2=-1$ and set $d_n=i r_{n+2}$.
Since $r_{n+2}^2=\lambda_{n+2}-1$, we have
$d_n^2=1-\lambda_{n+2}\ne0$. The maps
\[
 h_n:E_{v_{n+2}}\longrightarrow E_{T(v_n)},\qquad
 h_n(x,y)=\left(\frac{x-\lambda_{n+2}}{1-\lambda_{n+2}},
                    \frac{y}{d_n^3}\right)
\]
are isomorphisms over $\bar k$. They send $(\lambda_{n+2},0)$,
the generator of $\ker\widehat\Psi_{v_{n+1}}$, to $(0,0)$,
the generator of $\ker\Psi_{T(v_n)}$.

In fact, these choices give the exact compatibility
\[
 \Psi_{T(v_n)}\circ h_n
 =h_{n-1}\circ\widehat\Psi_{v_{n+1}}.
\]
Both sides are $2$-isogenies from $E_{v_{n+2}}$ to
$E_{T(v_{n-1})}$ with the same kernel. To check their
normalization, note that
$h_n^*\xi_{T(v_n)}=d_n\xi_{v_{n+2}}$ and that the dual of
$\Psi_{v_{n+1}}$ has differential multiplier $1/v_{n+1}$,
because its composition with $\Psi_{v_{n+1}}$ is $[2]$.
The two sides therefore have multipliers $2T(v_n)d_n$ and
$d_{n-1}/v_{n+1}$, respectively. They agree, since
\[
 2v_{n+1}T(v_n)\frac{r_{n+2}}{r_{n+1}}
 =T(v_n)\frac{(v_n-1)^2}{v_n^2-1}=1.
\]
Here we used $v_{n+1}^2-1=(v_n-1)^2/(2v_n)$.
Two isogenies with the same kernel differ by an automorphism
of their target fixing the origin; in characteristic $p\geq5$,
such an automorphism is the identity if its differential
multiplier is one. This proves the compatibility, including
the sign of the $y$-coordinate.

The reverse cycle as ordered in Proposition~\ref{prop:reverseroot}
starts at $T(v_{N-1})$, which is identified by $h_{N-1}$ with
$E_{v_1}$, rather than with $E_{v_0}$. To compare its composite
with $\widehat\Theta\in\End(E_{v_0})$, cyclically rotate the
reverse cycle so that it starts at $T(v_{N-2})$. If
$\Theta_{\mathrm{rev},0}$ denotes the resulting composite, the
compatibilities above give
\[
 \Theta_{\mathrm{rev},0}
 =h_{N-2}\circ\widehat\Theta\circ h_{N-2}^{-1}.
\]
Thus the appropriately based reverse composite is conjugate
to the dual by an isomorphism, with no remaining automorphism
ambiguity. Cyclic rotation does not change the product, and
the differential multiplier of $\widehat\Theta$ is
$t-2^N\omega=\omega^{-1}=2^N\omega^{\mathrm{rev},T}$,
in agreement with Proposition~\ref{prop:reverseroot}.
\end{remark}
\endgroup

%%%%%%%%%%%%%%%%%%%%%%%%%%%%%
\paragraph*{Reversal and Frobenius conjugation.}
%%%%%%%%%%%%%%%%%%%%%%%%%%%%%
The reverse cycle used above is obtained by reversing the order
\emph{and} applying $T$; reversing the order alone need not
preserve the recurrence. Its product is always $1/(2^N\omega)$.
By contrast, the Frobenius-conjugate cycle
$\mathbf v^{(p)}=(v_0^p,\ldots,v_{N-1}^p)$ has product $\omega^p$.
Since the product--trace quadratic has coefficients in $\FF_p$,
Frobenius permutes its roots. If $\omega\notin\FF_p$, it therefore
exchanges the two roots, giving
\[
 \omega^p=\frac{1}{2^N\omega}.
\]
If $\omega\in\FF_p$, Frobenius fixes this product, even when
the other root is distinct. Thus Frobenius conjugation and
$T$-reversal need not have the same effect on products.
The two operations commute, because $T(v^p)=T(v)^p$;
applying both gives a cycle with product $1/(2^N\omega^p)$.
This equals $\omega$ when $\omega\notin\FF_p$, and equals
$1/(2^N\omega)$ when $\omega\in\FF_p$.
These are statements about products: equality of products
does not assert equality of the corresponding cycles.

%%%%%%%%%%%%%%%%%%%%%%%%%%%%%
\section*{Acknowledgements}
%%%%%%%%%%%%%%%%%%%%%%%%%%%%%
This research has been conducted as part of GMV's Project
``Soluci\'on de Identidad Autosoberana'', and carried out under an
Art\'{\i}culo 60 contract at Universidad Carlos III de Madrid, with
partial co-financing from the Spanish National Cybersecurity Institute
(INCIBE) using resources from the Spanish ``Plan de Recuperaci\'on,
Transformaci\'on y Resiliencia'', originating from the European Union's
Recovery and Resilience Facility. The opinions expressed in this
publication are those of the authors and do not necessarily reflect
those of the co-financing entity.

Our warmest thanks to Mar\'{\i}a Isabel Gonz\'alez-Vasco, who promoted
this work, brought the authors together as principal investigator of
the UC3M node of the project, and read the manuscript and provided
valuable comments.

\end{document}